\documentclass[12pt]{article}
\usepackage{amsfonts}
\usepackage{full page,
amssymb, amscd, amsmath, graphicx, 
}
\usepackage{xcolor}
\newtheorem{thm}{Theorem}[section]
\newtheorem{defn}[thm]{Definition}
\newtheorem{prop}[thm]{Proposition}
\newtheorem{cor}[thm]{Corollary}
\newtheorem{lemma}[thm]{Lemma}
\newtheorem{rema}[thm]{Remark}

\newtheorem{exam}[thm]{Example}

\newcommand{\halmos}{\rule{1ex}{1.4ex}}

\newcommand{\nn}{\nonumber \\}

 \newcommand{\res}{\mbox{\rm Res}}

 \newcommand{\pf}{{\it Proof.}\hspace{2ex}}
 
 \newcommand{\epfv}{\hspace*{\fill}\mbox{$\halmos$}\vspace{1em}}

\newcommand{\lbar}{\bigg\vert}

\newcommand{\C}{\mathbb{C}}
\newcommand{\Z}{\mathbb{Z}}
\newcommand{\R}{\mathbb{R}}

\newcommand{\one}{\mathbf{1}}

\renewcommand{\mod}{\;\;\mbox{\rm mod}\;}
\newcommand{\g}{\mathfrak{g}}

\title{ {\bf  Weight-one elements of vertex operator algebras 
and automorphisms of categories of generalized 
twisted modules} }
\date{}
\author{Yi-Zhi Huang and Christopher Sadowski}

\begin{document}

\bibliographystyle{alpha}
\maketitle
\begin{abstract}
Given a weight-one element $u$ of a vertex operator algebra $V$, 
we construct an automorphism of the category of generalized $g$-twisted 
modules for automorphisms $g$ of $V$ fixing $u$.  
We apply this construction to the case that $V$ is an affine 
vertex operator algebra to obtain
explicit results on these automorphisms of categories. 
In particular, we give explicit constructions 
of certain generalized twisted modules from 
generalized twisted modules associated to diagram automorphisms
of finite-dimensional simple Lie algebras and generalized (untwisted) modules.  
\end{abstract}

\renewcommand{\theequation}{\thesection.\arabic{equation}}
\renewcommand{\thethm}{\thesection.\arabic{thm}}
\setcounter{equation}{0}
\setcounter{thm}{0}
\section{Introduction}

Orbifold conformal field theories are 
conformal field theories obtained from known conformal field theories
and their automorphisms. One conjecture in orbifold conformal 
field theory is that a category of suitable generalized twisted modules 
for a suitable vertex operator algebra has a natural structure of a
crossed tensor category in the sense of Turaev \cite{T}, where it is called 
a crossed group-category. See 
\cite{H-problems} and \cite{H-orbifold} 
for the precise conjectures and open problems. 

To study orbifold conformal field theories and to construct 
crossed tensor category structures, we first need to study categories of 
generalized twisted modules satisfying suitable conditions. In 
\cite{H-const-twisted-mod} and \cite{H-exist-twisted-mod}, lower-bounded 
generalized twisted modules for a grading-restricted vertex algebra
are constructed and studied. In \cite{H-aff-va-twisted-mod}, 
the constructions and results in \cite{H-const-twisted-mod} 
and \cite{H-exist-twisted-mod} are applied to affine vertex (operator) algebras
to give explicit constructions of various types 
of lower-bounded and grading-restricted generalized 
twisted modules.

In this paper, for each weight-one element $u$ of a 
vertex operator algebra $V$, 
we first construct a transformation
sending a generalized $g$-twisted 
module for an automorphism $g$ of $V$ fixing $u$ to another
such generalized twisted module by generalizing the construction 
of examples of generalized twisted 
modules in \cite{H-log-twisted-mod}. In the special case that the automorphisms
of $V$ involved (including the automorphisms obtained from 
the zero modes of the vertex operators
of weight-one elements) are of finite order, 
such transformations for weak twisted modules (without gradings) 
were constructed by Li \cite{L}. We then show that such a transformation
is in fact an automorphism of the category of generalized $g$-twisted 
modules for automorphisms $g$ fixing $u$. 
Such automorphisms of categories 
seem to be independent of the conjectured crossed tensor category structure of 
suitable generalized twisted modules. 
These automorphisms should be a special feature of the conjectured 
crossed tensor category structure of suitable generalized twisted 
$V$-modules and are expected to be a necessary part  
in the tensor categorical study of 
orbifold conformal field theories, especially in the study of problems
related to ``crossed vertex tensor categories'' (the twisted generalizations of 
vertex tensor categories introduced in \cite{HL}) and the reconstruction
problem. 

We apply our construction mentioned above to the case of affine 
vertex operator algebras to obtain
explicit results on these automorphisms of categories. 
In particular, we give explicit constructions 
of certain generalized twisted modules from 
generalized twisted modules associated to diagram automorphisms
of finite-dimensional simple Lie algebras and generalized (untwisted) modules. 

This paper is organized as follows: In Section 2,
we review
some basic facts on automorphisms 
of $V$ and recall a formal series $\Delta_{V}^{(u)}$  of operators 
on $V$ associated to 
a weight-one element $u$ of $V$ introduced in \cite{H-log-twisted-mod}.  
In Section 3, we recall 
the notions of variants of generalized twisted module for a 
vertex operator algebra $V$. 
We construct in Section 4 an automorphism of
the category of generalized $g$-twisted $V$-modules for all automorphisms 
$g$ of $V$ fixing $u$. We recall some basic facts about affine vertex 
(operator) algebras $M(\ell, 0)$ and $L(\ell, 0)$
and their automorphisms in Section 5. 
In Sections 6 and 7, we discuss the applications of the 
automorphisms of the categories of generalized twisted modules 
to the case of affine vertex operator algebras $M(\ell, 0)$ and $L(\ell, 0)$. 
In Section 6, we construct explicitly 
generalized $g$-twisted $M(\ell, 0)$- and $L(\ell, 0)$-modules 
for  semisimple
automorphisms $g$ of $M(\ell, 0)$- and $L(\ell, 0)$ 
from generalized $\mu$-twisted modules for diagram automorphisms 
$\mu$ of the finite-dimensional simple Lie algebra. 
In particular, we construct 
generalized $g$-twisted modules for semisimple 
automorphisms $g$ of $M(\ell, 0)$ and $L(\ell, 0)$ 
from generalized (untwisted) $V$-modules. 
In Section 7, we construct explicitly 
generalized $g$-twisted $M(\ell, 0)$- and $L(\ell, 0)$-modules 
for general (not-necessarily semisimple)
automorphisms $g$ of $M(\ell, 0)$ and $L(\ell, 0)$ 
from generalized $\mu$-twisted modules for diagram automorphisms 
$\mu$ of the finite-dimensional simple Lie algebra. 
In particular, we construct 
generalized $g$-twisted modules for general
automorphisms $g$ of $M(\ell, 0)$ and $L(\ell, 0)$ 
from generalized (untwisted) $V$-modules. 

\renewcommand{\theequation}{\thesection.\arabic{equation}}
\renewcommand{\thethm}{\thesection.\arabic{thm}}
\setcounter{equation}{0}
\setcounter{thm}{0}

\section{Automorphisms of $V$ and formal 
series $\Delta_{V}^{(u)}$ of operators}

In this section, we review
some basic facts about automorphisms 
of a vertex operator algebra $V$
and recall a formal series $\Delta_{V}^{(u)}$  of operators 
on $V$ associated to 
a weight-one element $u$ of $V$ introduced in \cite{H-log-twisted-mod}. 

Let $g$ be an automorphism of $V$, that is, 
a linear isomorphism $g: V\to V$ such that 
$g(Y_{V}(u, x)v)=Y_{V}(g(u), x)g(v)$, $g(\one)=\one$ and 
$g(\omega)=\omega$. 
Using the multiplicative Jordan-Chevalley decomposition, 
there are unique commuting operators $\sigma$ and 
$g_\mathcal{U}$ such that $g=\sigma g_\mathcal{U}$, where $\sigma$ is a 
semisimple automorphism of $V$ and $g_\mathcal{U}$ is a locally unipotent 
operator on $V$, in the sense that $g_\mathcal{U}-1_V$ is 
locally nilpotent on $V$.
We have that $V = \coprod_{\alpha \in P_{V}^{g}} V^{[\alpha]}$, where
\begin{align*}
P_{V}^{g} = \{ \alpha \in [0,1) + i\R \mid e^{2 \pi i \alpha} 
\text{ is an eigenvalue of } g \}.
\end{align*}
Following \cite{HY}, we write $\sigma = e^{2 \pi i \mathcal{S}_g}$, where 
$\mathcal{S}_g$ is the operator defined by $\mathcal{S}_g v = \alpha v$ 
for $v \in V^{[\alpha]}$ for each $\alpha \in P_{V}^{g}$ and extended 
linearly, and we write $g_\mathcal{U} = e^{2 \pi i \mathcal{N}_g}$, 
where $$\mathcal{N}_g = \sum_{j \in \mathbb{Z}_+} \frac{(-1)^j 
(g_\mathcal{U}-1_V)^j}{2 \pi i j}.$$ We note that the generalized 
eigenspaces of $g$ are the eigenspaces $V^{[\alpha]}$ 
of $\sigma$ and $\mathcal{S}_g$ for 
$\alpha\in P_{V}^{g}$, and that $\mathcal{N}_g$ is locally nilpotent since 
$g_\mathcal{U}-1_V$ is locally nilpotent. Also, we note that since $\sigma$ 
and $g_\mathcal{U}$ commute, they preserve 
each other's generalized eigenspaces. In particular, if $v \in V$ is an eigenvector of $\sigma$ with 
eigenvalue $e^{2 \pi i \alpha}$, then so is $g_\mathcal{U} v$. Since the 
eigenspaces of $\sigma$ are the same as the eigenspaces of $\mathcal{S}_g$, 
we have that if $v$ is an eigenvector of $\mathcal{S}_g$ with eigenvalue 
$\alpha$, then so is $g_{\mathcal{U}} v$. Hence, by the definition of 
$\mathcal{N}_g$, we see that $\mathcal{N}_g v$ is also an eigenvector of 
$\mathcal{S}_g$ with eigenvalue $\alpha$. Thus, on an eigenvector $v$ 
of $\mathcal{S}_g$ with eigenvalue $\alpha$, 
\begin{align*}
\mathcal{S}_g \mathcal{N}_g v = \alpha \mathcal{N}_g v 
= \mathcal{N}_g (\alpha v) = \mathcal{N}_g \mathcal{S}_g v
\end{align*}
and so $[\mathcal{S}_g,\mathcal{N}_g] = 0$ on $V$. 
Hence, we have $g = e^{2 \pi i \mathcal{L}_g}$ where 
$\mathcal{L}_g = \mathcal{S}_g + \mathcal{N}_g$.

For $u \in V$, 
we write $Y(u,x) = \sum_{n \in \mathbb{Z}} Y_n(u) x^{-n-1}$.
Recall from Proposition 5.1 in \cite{H-log-twisted-mod}, 
for $u \in V_{(1)}$  satisfying $L(1)u=0$, 
$g_u = e^{2 \pi i Y_0(u)}$ is an automorphism of the vertex operator algebra
$V$.
Using the notation in 
\cite{H-const-twisted-mod}, we may decompose $g_u$ into its
semisimple and unipotent parts, which we denote $e^{2 \pi i Y_0(u)_{S}}$
and $e^{2 \pi i Y_0(u)_{N}}$, respectively. Here $Y_0(u)_{S}$ is the semisimple 
part of $Y_0(u)$ and $Y_0(u)_{N}$ is the nilpotent part of $Y_0(u)$ so that
$$Y_0(u) = Y_0(u)_S + Y_0(u)_N.$$

We also recall from \cite{H-log-twisted-mod} 
the formal series $\Delta^{(u)}_V(x)$ for $u\in V_{(1)}$ defined by
$$\Delta_{V}^{(u)}(x)=x^{-Y_{0}(u)}e^{-\int_{0}^{-x}Y^{\le -2}(u, y)} 
= x^{-Y_0(u)_{S}} e^{-2 \pi i (Y_0(u)_{N}) \log x} 
e^{-\int_{0}^{-x}Y^{\le -2}(u, y)} $$
where the linear map 
$$\int_{0}^{x}: V[[x]]+x^{-2}V[[x^{-1}]]
\to xV[[x]]+x^{-1}V[[x^{-1}]]$$
is defined by 
$$\int_{0}^{x} \sum_{n\in \Z\setminus \{-1\}}v_{n}y^{n}
=\sum_{n\in \Z\setminus \{-1\}}\frac{v_{n}}{n+1}x^{n+1}$$
and where $Y^{\le -2}(u,x) = \sum_{m \in \mathbb{Z}_+} Y_m(u) x^{-m-1}$. 
The properties of $\Delta^{(u)}_V(x)$ for $u\in V_{(1)}$ in
the following theorem are given in Theorem 5.2 in \cite{H-log-twisted-mod}
except for (\ref{D3}):

\begin{thm}\label{H1-Delta-Theorem}
For $v \in V$, there exist $m_1,\dots,m_l \in \mathbb{R}$ such that 
\begin{align}\label{D1}
\Delta^{(u)}_V(x)v \in x^{m_1}V[x^{-1}][\log x] + \cdots + x^{m_l}V[x^{-1}][\log x].
\end{align}
Moreover, the series $\Delta^{(u)}_V(x)$ satisfies
\begin{align}
\Delta^{(u)}_V(x)Y(v,x_2) &= Y(\Delta^{(u)}_V(x+x_2)v,x_2)\Delta^{(u)}_V(x),
\label{D2}\\
[L(0),\Delta^{(u)}_V(x)] &= x \frac{d}{dx}\Delta^{(u)}_V(x) + Y_0(u)\Delta^{(u)}_V(x), \label{D3}\\
[L(-1),\Delta^{(u)}_V(x)] &= -\frac{d}{dx}\Delta^{(u)}_V(x).\label{D4}
\end{align}
\end{thm}
\pf
We need only prove (\ref{D3}). 
Using the fact that $u \in V_{(1)}$, we have by the commutator formula that
\begin{equation}
[L(0),Y_m(u)] = -mY_m(u)
\end{equation}
for $u \in V_{(1)}$. Thus, we have that, for $v \in V$, 
\begin{align*}
[L(0),\Delta^{(u)}_V(x)] &= \left[L(0),x^{-Y_0(u)}\exp \left(-\sum_{m \ge 1} 
\frac{Y_m(u)}{-m}(-x)^{-m}\right)\right]\\
&=x^{-Y_0(u)} \left[ L(0), \sum_{k \ge 0} \frac{1}{k!} \left( -\sum_{m \ge 0}  
\frac{(-1)^m Y_m(u)}{-m}x^{-m}  \right)^k \right]\\
&= x^{-Y_0(u)}\left[ L(0), \sum_{k \ge 0} \frac{1}{k!} 
\left( -\sum_{\substack{n_1,\dots,n_k \in \Z_+\\ n_1 + \cdots + n_k = l}}  
\frac{(-1)^{l} Y_{n_1}(u)\cdots Y_{n_k}(u)}{(-n_1)\cdots(-n_k)} x^{-l}  
\right) \right]\\
&=  x^{-Y_0(u)}\sum_{k \ge 0} \frac{1}{k!} 
\left( -\sum_{\substack{n_1,\dots,n_k \in \Z_+\\ n_1 + \cdots + n_k = l}}
(-n_1 - \cdots - n_k)  \frac{(-1)^{l} Y_{n_1}(u)\cdots Y_{n_k}(u)}{(-n_1)
\cdots(-n_k)} x^{-l}  \right)\\
&=x^{-Y_0(u)}x\frac{d}{dx}\sum_{k \ge 0} \frac{1}{k!} 
\left( -\sum_{\substack{n_1,\dots,n_k \in \Z_+\\ n_1 + \cdots + n_k = l}} 
\frac{(-1)^{l} Y_{n_1}(u)\cdots Y_{n_k}(u)}{(-n_1)\cdots(-n_k)} x^{-l}  
\right)\\
&= x^{-Y_0(u)} x\frac{d}{dx}\exp \left( -\sum_{m \ge 1} 
\frac{Y_m(u)}{-m}(-x)^{-m}\right)\\
&= x\frac{d}{dx}\Delta^{(u)}_V(x) + Y_0(u)\Delta^{(u)}_V(x).
\end{align*}
\epfv

We also note the following properties of $\Delta^{(u)}_V(x)$:

\begin{prop}\label{Delta-prop}
Let $\Delta^{(u)}_V(x)$ for $u\in V_{(1)}$ be the series defined above.
\begin{enumerate}
\item $\Delta^{(0)}(x) = 1_V$
\item If $u_1, u_2 \in V$ and $[Y(u_1,x_1),Y(u_2,x_2)]=0$, then
$$ \Delta^{(u_1)}_V(x)\Delta^{(u_2)}_V(x) = \Delta^{(u_1 +u_2)}_V(x)$$
\item $\Delta^{(u)}_V(x)$ is invertible and $\Delta^{(u)}_V(x)^{-1} 
= \Delta^{(-u)}_V(x)$.
\end{enumerate}
\end{prop}
\pf
Property 1 is obvious. Property 2 follow from properties of the 
exponential function. Property 3 is a consequence of Properties 1 and 2.
\epfv

In the proof of the main theorem in Section 4,
we also need a lemma.

\begin{lemma}
Let $g$ be an automorphism of $V$ and suppose $u \in V_{(1)}$ 
such that $g(u) = u$. Then, we have that 
$$
[g,\Delta^{(u)}_V(x)] = 0.
$$
\end{lemma}
\pf
First, we note that since $g$ is an automorphism of $V$, we have that
$$
gY(u,x) = Y(g(u),x)g = Y(u,x)g
$$
so that 
$$
[g,Y(u,x)] = 0
$$
and so $[g,Y_n(u)]=0$ for all $n \in \Z$.
In particular, we have that $[g,Y_0(u)]=0$ and so $[g,g_u] = 0$. Next, 
we show that $$[g, Y_0(u)_{S}] = [g,Y_0(u)_N] = 0.$$
Since $[g,Y_0(u)]=0$, we have that $g(Y_0(u)-\lambda I)^k 
= (Y_0(u)-\lambda I)^k g$ for any $k \in \mathbb{N}$ and  
$\lambda \in \mathbb{C}$, so that $g$ preserves the generalized 
eigenspaces of $Y_0(u)$. Thus, if $v$ is a generalized eigenvector 
of $Y_0(u)$ with eigenvalue $\lambda \in \mathbb{C}$, we have 
that $g Y_0(u)_{S} v = g\lambda v$ and $Y_0(u)_{S}gv = \lambda g v$ 
and so $[g,Y_0(u)_{S}] = 0$. From this we immediately have 
$[g,Y_0(u)_{N}] = 0$. Our lemma now follows.
\epfv

Using the fact that $[g,Y_n(u)]=0$ for each $n \in \mathbb{Z}$, 
a similar proof also shows that $Y_n(u)$ preserves the generalized 
eigenspaces $V^{[\alpha]}$ of $g$ for $\alpha \in P_{V}^{g}$. 
Hence, $Y_n(u)$ preserves the eigenspaces of $\mathcal{S}_g$ 
and $e^{2 \pi \mathcal{S}_g}$, and so we have that 
$[\mathcal{S}_g,Y_n(u)] =[e^{2 \pi i \mathcal{S}_g},Y_n(u)]= 0$ 
on $V^{[\alpha]}$ and hence on $V$. From this, it follows that
$[e^{-2 \pi i \mathcal{S}_g} g, Y_n(u)] = [g_\mathcal{U},Y_n(u)] = 0$.
Using the fact that $\mathcal{N}_g = \sum_{j \in \mathbb{Z}_+} 
\frac{(-1)^j (g_\mathcal{U}-1_V)^j}{2 \pi i j}$, we have 
$[\mathcal{N}_g,Y_n(u)] = 0$. Hence, we have 
$[\mathcal{L}_g, Y_n(u)] =0$
for all $n \in \mathbb{Z}$.


\renewcommand{\theequation}{\thesection.\arabic{equation}}
\renewcommand{\thethm}{\thesection.\arabic{thm}}
\setcounter{equation}{0}
\setcounter{thm}{0}

\section{Generalized twisted modules}

We recall the notions of variants of generalized twisted module for a 
vertex operator algebra in this section. Our definition 
uses the Jacobi identity for the twisted vertex operators formulated in 
\cite{H-twist-vo} as the main axiom. See \cite{H-log-twisted-mod} for 
the original definition in terms of duality properties, \cite{B} for the 
component form of the Jacobi identity and \cite{HY} for a derivation
of the Jacobi identity from the duality properties.

\begin{defn}\label{C/Z-g-twisted-mod}
{\rm Let $V$ be a vertex operator algebra and 
$g$ an automorphism of $V$.
A {\it $\mathbb{C}/\mathbb{Z}$-graded generalized 
$g$-twisted $V$-module} is a 
$\mathbb{C} \times \mathbb{C}/\mathbb{Z}$-graded vector space 
$W = \coprod_{n \in \mathbb{C}, \alpha \in \mathbb{C}/\mathbb{Z}} 
W_{[n]}^{[\alpha]}$ equipped with an action of $g$
and a linear map
\begin{align*}
Y^g_W: V\otimes W &\to W\{x\}[{\rm log} x],\\
v \otimes w &\mapsto Y^g_W(v, x)w
\end{align*}
satisfying:
\begin{enumerate}
\item The {\it equivariance property}: For $p \in \mathbb{Z}$, 
$z \in \mathbb{C}^\times$, $v \in V$ and $w \in W$, we have
$$
Y^{g;p+1}_W(gv,z)w = Y^{g;p}_W(v,z)w
$$
where 
$$Y^{g; p}_W(v, z)w=Y^{g}_W(v, x)w\lbar_{x^{n}=e^{nl_{p}(z)},
\; \log x=l_{p}(z)}$$
is the $p$-th analytic branch of $Y^g$.

\item The {\it identity property}: For $w \in W$, $Y^g_W({\bf 1}, x)w
= w$.

\item The {\it lower truncation property}: For all $v \in V$ and $w \in W$, 
$Y^g_W(v,x)w$ is lower truncated, that is, $(Y^g_W)_n(v)w = 0$ 
when $\Re(n)$ is sufficiently negative.

\item The {\it Jacobi identity}: \begin{eqnarray}\label{Jacobi}
x_0^{-1}\delta\left(\frac{x_1 - x_2}{x_0}\right)Y^{g}_W(u, x_1)Y^{g}_W(v, x_2)
- x_0^{-1}\delta\left(\frac{- x_2 + x_1}{x_0}\right)Y^{g}_W(v, x_2)Y^{g}_W(u, x_1)\nn
= x_1^{-1}\delta\left(\frac{x_2+x_0}{x_1}\right)
Y^{g}_W\left(Y\left(\left(\frac{x_2+x_0}{x_1}\right)^{{\mathcal{L}}_g}u, x_0\right)v, x_2\right)
\end{eqnarray}

\item Properties about the gradings: Let $L^{g}_{W}(0)=\res_{x}xY^{g}_{W}(\omega, x)$.
Then we have:
\begin{enumerate}
\item The {\it $L(0)$-grading condition}: For each $w \in W_{[n]} 
= \coprod_{\alpha \in \C / \Z} W^{[\alpha]}_{[n]}$,
there exists $K \in \mathbb{Z}_+$ such that 
\begin{equation*}
(L^g_W(0) - n)^K w = 0.
\end{equation*}
\item The {\it $g$-grading condition}: For each 
$w \in W^{[\alpha]} = \coprod_{n \in \C} W^{[\alpha]}_{[n]}$,
there exists $\Lambda \in \mathbb{Z}_+$ such that
\begin{equation*}
(g - e^{2 \pi i \alpha})^\Lambda w = 0
\end{equation*}
\end{enumerate}
where $Y_W^g(\omega,x) = \sum_{n \in \mathbb{Z}} L^g_W(n) x^{-n-2}$.
\item The $L(-1)$-{\it derivative property}: 
For $v \in V$,
$$
\frac{d}{dx}Y^g_W(v, x) = Y^g_W(L_{V}(-1)v, x).
$$
\end{enumerate}
}
\end{defn}

We denote the $\mathbb{C}/\mathbb{Z}$-graded generalized 
$g$-twisted $V$-module defined above by
$(W, Y^{g}_{W})$) or simply by $W$. 

\begin{defn}
\rm{We call a generalized $g$-twisted $V$-module {\it lower-bounded} 
if  $W_{[n]} = 0$ for $n\in \C$ such that 
$\Re(n)$ is sufficiently negative. We call 
a lower bounded generalized $g$-twisted $V$-module {\it grading-restricted} 
(or simply a {\it $g$-twisted $V$-module}) if for each $n \in \C$ we have 
$\text{dim} W_{[n]} < \infty$.}
\end{defn}

Throughout the work, we will write
$$W = \coprod_{n \in \mathbb{C}, \alpha \in P^g_W} W_{[n]}^{[\alpha]} $$
where $$ P^g_W = \{\Re(\alpha) \in [0,1) | e^{2 \pi i \alpha} 
\text{ is an eigenvalue of } g \text{ acting on } W \}.$$

\begin{defn}
{\rm A vertex operator algebra $V$ 
is said to have 
a {\it (strongly)
$\mathbb{C}$-graded vertex operator algebra structure compatible with $g$} 
if $V$ has an additional $\C$-grading
$V=\coprod_{\alpha\in \C}V^{[\alpha]}
=\coprod_{n\in \Z,\; \alpha\in \C}V_{(n)}^{[\alpha]}$ such that 
for $V^{[\alpha]}$ for $\alpha\in \C$ is the generalized eigenspace of $g$ with 
eigenvalue $e^{2\pi i\alpha}$. For such a vertex operator algebra, 
a {\it $\mathbb{C}$-graded generalized $g$-twisted $V$-module} is a 
$\mathbb{C} \times \mathbb{C}$-graded vector space 
$W = \coprod_{n,\alpha \in \mathbb{C}} W_{[n]}^{[\alpha]}$ 
equipped with an action of $g$ and a linear map
\begin{align*}
Y^g_W: V\otimes W &\to W\{x\}[{\rm log} x],\\
v \otimes w &\mapsto Y^g_W(v, x)w
\end{align*}
satisfying the same axioms as in Definition \ref{C/Z-g-twisted-mod} except 
that $\mathbb{C}/\mathbb{Z}$ is replaced by $\mathbb{C}$ 
and the {\it grading compatibility condition} holds: For 
$\alpha, \beta \in \mathbb{C}, v \in V^{[\alpha]},$ 
and $w \in W^{[\beta]}$ we have
$$ Y^g_W(v,x)w \in W^{[\alpha + \beta]}[\log x].$$
}
\end{defn}

\begin{defn}
{\rm
A generalized $\mathbb{C}/\mathbb{Z}$-graded (or $\mathbb{C}$-graded) 
$g$-twisted $V$-module $W$ is said to be  {\it strongly 
$\mathbb{C}/\mathbb{Z}$-graded} (or {\it strongly $\mathbb{C}$-graded}) 
if it satisfies the {\it grading restriction condition}: For each $n \in \mathbb{C}$ 
and $\alpha \in \mathbb{C}/\mathbb{Z}$ (or $\alpha \in \mathbb{C}$) 
dim$W_{[n]}^{[\alpha]} < \infty$ and $W^{[\alpha]}_{[n +l ]} = 0$ 
for all sufficiently negative real numbers $l$.}
\end{defn}

The original definition of lower-bounded 
generalized $V$-module in \cite{H-log-twisted-mod} is in terms of the duality property. 
More precisely, for a lower-bounded generalized $g$-twisted $V$-module
$W$, it was shown in \cite{HY} that 
the lower truncation property for the twisted vertex operator map
and Jacobi identity is equivalent to the duality property:
 Let $W^{'} = \coprod_{n \in \mathbb{C}, \alpha \in
\mathbb{C}/\mathbb{Z}} \left(W_{[n]}^{[\alpha]}\right)^{\ast}$ and, for $n \in
\mathbb{C}$, $\pi_n : W \rightarrow W_{[n]}=\coprod_{\alpha\in \C/\Z}
W_{[n]}^{[\alpha]}$ be the projection. For
any $u, v \in V$, $w \in W$ and $w^{'} \in W^{'}$, there exists a
multivalued analytic function of the form
\[
f(z_1, z_2) = \sum_{i,
j, k, l = 0}^N a_{ijkl}z_1^{m_i}z_2^{n_j}({\rm log}z_1)^k({\rm
log}z_2)^l(z_1 - z_2)^{-t}
\]
for $N \in \mathbb{N}$, $m_1, \dots,
m_N$, $n_1, \dots, n_N \in \mathbb{C}$ and $t \in \mathbb{Z}_{+}$,
such that the series
\[
\langle w^{'}, Y^{g; p}(u, z_1)Y^{g;p}(v,
z_2)w\rangle = \sum_{n \in \mathbb{C}}\langle w^{'}, Y^{g; p}(u,
z_1)\pi_nY^{g; p}(v, z_2)w\rangle,
\]
\[
\langle w^{'}, Y^{g; p}(v,
z_2)Y^{g;p}(u, z_1)w\rangle = \sum_{n \in \mathbb{C}}\langle w^{'},
Y^{g; p}(v, z_2)\pi_nY^{g; p}(u, z_1)w\rangle,
\]
\[
\langle w^{'}, Y^{g; p}(Y(u, z_1 - z_2)v,
z_2)w\rangle = \sum_{n \in \mathbb{C}}\langle w^{'}, Y^{g;
p}(\pi_nY(u, z_1 - z_2)v, z_2)w\rangle
\]
are absolutely convergent in
the regions $|z_1| > |z_2| > 0$, $|z_2| > |z_1| > 0$, $|z_2| > |z_1
- z_2| > 0$, respectively, and are convergent to the branch
\[
 \sum_{i, j, k, l = 0}^N
a_{ijkl}e^{m_il_p(z_1)}e^{n_jl_p(z_2)}l_p(z_1)^kl_p(z_2)^l(z_1 -
z_2)^{-t}
\]
of $f(z_1, z_2)$ when $\arg z_{1}$ and $\arg z_{2}$ are sufficiently close (more precisely, 
when $|\arg z_{1}-\arg z_{2}|<\frac{\pi}{2}$).

\renewcommand{\theequation}{\thesection.\arabic{equation}}
\renewcommand{\thethm}{\thesection.\arabic{thm}}
\setcounter{equation}{0}
\setcounter{thm}{0}

\section{Automorphisms of categories of generalized twisted modules}

Let $V$ be a vertex operator algebra. 
In this section, given $u\in V_{(1)}$,
we generalize the construction of generalized twisted 
$V$-modules from (untwisted) generalized $V$-modules in 
\cite{H-log-twisted-mod} to obtain a transformation which sends 
generalized $g$-twisted $V$-modules for automorphisms $g$ of 
$V$ fixing $u$ to 
generalized twisted $V$-modules of the same type. 
We then show that this transformation is in fact an automorphism of
the category of generalized $g$-twisted $V$-modules for
automorphisms $g$ of $V$ fixing $u$.

We first prove a result on the properties of a modified vertex operator map
from a $\C/\Z$-graded (or $\C$-graded) generalized 
$g$-twisted $V$-module using $\Delta_{V}^{(u)}(x)$.

\begin{thm}\label{main-thm}
Let $V$ be a vertex operator algebra satisfying 
$V_{(0)}=\mathbb{C}\one$, $V_{(n)} = 0$ 
for $n < 0$ and $L(1)u = 0$, $g$ an automorphism of $V$ and 
$u \in V_{(1)}$ satisfying 
$g(u) = u$.
Let $(W, Y_{W}^{g})$ 
be a $\C/\Z$-graded (or $\C$-graded) generalized 
$g$-twisted $V$-module. Then the map 
\begin{align*}
Y_W^{g g_u}: V\otimes W&\to W\{x\}[\log x]\\
v\otimes w&\mapsto Y_{W}^{g}(v, x)w
\end{align*} 
defined by
$$
Y^{gg_u}_W(v,x) = Y^g_W(\Delta^{(u)}_V(x) v, x)
$$
for $v\in V$ satisfies the identity property, the lower truncation property, 
the $L(-1)$-derivative property, the equivariance property
and the Jacobi identity. 
\end{thm}
\pf
The proof of the identity property and $L(-1)$-derivative property are 
identical to the proof in \cite{H-log-twisted-mod}. The lower truncation 
property of 
$Y^{gg_u}_W$ follows from Theorem \ref{H1-Delta-Theorem} 
and the lower truncation property of $Y^g_W$. 

We now prove that the equivariance property holds, that is, we show that 
$$Y^{gg_u;p+1}(gg_uv,z)w= Y^{gg_u;p}(v,z)w$$
for $v \in V$ and $w \in W$. Recall from the preceding section that since $g(u) = u$, we have that 
$[g,Y_n(u)] = 0$ for all $n \in \mathbb{Z}$. Let $v \in V$ be a generalized eigenvector 
of $Y_0(u)$ with eigenvalue $\lambda$. Then, we have
\begin{align*}
\Delta&^{(u)}_V(x) g g_u v\lbar_{x^{n}=e^{nl_{p+1}(z)},\; 
\log x=l_{p+1}(z)}\nn
 &= g \Delta^{(u)}_V(x) e^{2 \pi i Y_0(u)} v \lbar_{x^{n}
=e^{nl_{p+1}(z)},\; \log x=l_{p+1}(z)}\\
&= g e^{-\int_{0}^{-x}Y^{\le -2}(u, y)} x^{-\lambda}
e^{(-Y_0(u)_{N})\log x}e^{2 \pi i \lambda} e^{2 \pi i (Y_0(u)_N)}v 
\lbar_{x^{n}=e^{nl_{p+1}(z)},\; \log x=l_{p+1}(z)}\\
&= g e^{-\int_{0}^{-x}Y^{\le -2}(u, y)} 
 x^{-\lambda} e^{(-Y_0(u)_{N})\log x} v \lbar_{x^{n}
=e^{nl_{p}(z)},\; \log x=l_{p}(z)}\\
&=g \Delta^{(u)}_V(x)v\lbar_{x^{n}=e^{nl_{p}(z)},\; 
\log x=l_{p}(z)}.
\end{align*}
Thus, we have 
\begin{align*}
Y^{gg_u; p+1}_W(gg_uv,z)w &= Y^{g}_W
(\Delta^{(u)}_V(x)gg_uv,x)w \lbar_{x^{n}
=e^{nl_{p+1}(z)},\; \log x=l_{p+1}(z)}\\
&= Y^{g}_W\left(g\Delta^{(u)}_V(x)v \lbar_{x^{n}
=e^{nl_{p}(z)},\; \log x=l_{p}(z)},x\right)w \lbar_{x^{n}
=e^{nl_{p+1}(z)},\; \log x=l_{p+1}(z)}\\
&= Y^{g}_W(\Delta^{(u)}_V(x)v,x)w \lbar_{x^{n}
=e^{nl_{p}(z)},\; \log x=l_{p}(z)}\\
&=Y^{gg_u;p}(v,z)w.
\end{align*}

Finally, we show that the Jacobi identity holds. Using the 
Jacobi identity for $Y^g_W$, we have for, $v_1,v_2 \in V$,
\begin{align*}
&x_0^{-1}\delta\left(\frac{x_1 - x_2}{x_0}\right)Y^{gg_u}_W
(v_1, x_1)Y^{gg_u}_W(v_2, x_2)
- x_0^{-1}\delta\left(\frac{- x_2 + x_1}{x_0}\right)Y^{gg_u}_W
(v_2, x_2)Y^{gg_u}_W(v_1, x_1)\\
&=x_0^{-1}\delta\left(\frac{x_1 - x_2}{x_0}\right)Y^{g}_W
(\Delta^{(u)}_V(x_1)v_1, x_1)Y^{g}_W(\Delta^{(u)}_V(x_2)
v_2, x_2)\nn
&\quad - x_0^{-1}\delta\left(\frac{- x_2 + x_1}{x_0}\right)
Y^{g}_W(\Delta^{(u)}_V(x_2)v_2, x_2)Y^{g}_W
(\Delta^{(u)}_V(x_1)v_1, x_1)\nn
&= x_1^{-1}\delta\left(\frac{x_2+x_0}{x_1}\right)
Y^{g}_W\left(Y\left(\left(\frac{x_2+x_0}{x_1}\right)
^{{\mathcal{L}}_g}\Delta^{(u)}_V(x_1)v_1, x_0\right)
\Delta^{(u)}_V(x_2)v_2, x_2\right).
\end{align*}
Our desired Jacobi identity is:
\begin{align*}
&x_0^{-1}\delta\left(\frac{x_1 - x_2}{x_0}\right)Y^{gg_u}_W
(v_1, x_1)Y^{gg_u}_W(v_2, x_2)
- x_0^{-1}\delta\left(\frac{- x_2 + x_1}{x_0}\right)Y^{gg_u}_W
(v_2, x_2)Y^{gg_u}_W(v_1, x_1)\\
&= x_1^{-1}\delta\left(\frac{x_2+x_0}{x_1}\right)
Y^{gg_u}_W\left(Y\left(\left(\frac{x_2+x_0}{x_1}\right)
^{{\mathcal{L}}_{gg_u}}v_1, x_0\right)v_2, x_2\right)\\
&= x_1^{-1}\delta\left(\frac{x_2+x_0}{x_1}\right)
Y^{g}_W\left(\Delta^{(u)}_V(x_2)Y\left(\left(\frac{x_2+x_0}
{x_1}\right)^{{\mathcal{L}}_{g}+Y_0(u)}v_1, x_0\right)
v_2, x_2\right)\\
\end{align*}
Thus, it suffices to prove that
\begin{align*}
& x_1^{-1}\delta\left(\frac{x_2+x_0}{x_1}\right)
Y\left(\left(\frac{x_2+x_0}{x_1}\right)^{{\mathcal{L}}_g}
\Delta^{(u)}_V(x_1)v_1, x_0\right)\Delta^{(u)}_V(x_2)\\
&= x_1^{-1}\delta\left(\frac{x_2+x_0}{x_1}\right)
\Delta^{(u)}_V(x_2)Y\left(\left(\frac{x_2+x_0}{x_1}\right)
^{{\mathcal{L}}_{g}+Y_0(u)}v_1, x_0\right)\\
\end{align*}
or equivalently
\begin{align*}
& x_1^{-1}\delta\left(\frac{x_2+x_0}{x_1}\right)
Y\left(\left(\frac{x_2+x_0}{x_1}\right)^{{\mathcal{L}}_g}
\Delta^{(u)}_V(x_1)v_1, x_0\right)\\
&= x_1^{-1}\delta\left(\frac{x_2+x_0}{x_1}\right)
\Delta^{(u)}_V(x_2)Y\left(\left(\frac{x_2+x_0}{x_1}\right)
^{{\mathcal{L}}_{g}+Y_0(u)}v_1, x_0\right)\Delta^{(u)}_V
(x_2)^{-1}.\\
\end{align*}
By Theorem \ref{H1-Delta-Theorem} we have
\begin{align*}
&x_1^{-1}\delta\left(\frac{x_2+x_0}{x_1}\right)
\Delta^{(u)}_V(x_2)Y\left(\left(\frac{x_2+x_0}{x_1}\right)
^{{\mathcal{L}}_{g}+Y_0(u)}v_1, x_0\right)\Delta^{(u)}_V(x_2)^{-1}\\
&= x_1^{-1}\delta\left(\frac{x_2+x_0}{x_1}\right)
Y\left(\Delta^{(u)}_V(x_2+x_0)\left(\frac{x_2+x_0}{x_1}\right)
^{{\mathcal{L}}_{g}+Y_0(u)}v_1, x_0\right).
\end{align*}
So it suffices to prove
\begin{align*}
& x_1^{-1}\delta\left(\frac{x_2+x_0}{x_1}\right)
\left(\frac{x_2+x_0}{x_1}\right)^{{\mathcal{L}}_g}
\Delta^{(u)}_V(x_1)v_1\\
&= x_1^{-1}\delta\left(\frac{x_2+x_0}{x_1}\right)
\Delta^{(u)}_V(x_2+x_0)\left(\frac{x_2+x_0}{x_1}\right)
^{{\mathcal{L}}_{g}+Y_0(u)}v_1
\end{align*}
In fact,
\begin{align*}
 x_1^{-1}&\delta\left(\frac{x_2+x_0}{x_1}\right)
\Delta^{(u)}_V(x_2+x_0)\left(\frac{x_2+x_0}{x_1}\right)
^{{\mathcal{L}}_{g}+Y_0(u)}v_1\nn
&= x_1^{-1}\delta\left(\frac{x_2+x_0}{x_1}\right)
(x_2+x_0)^{-Y_0(u)} e^{-\int_{0}^{-(x_2+x_0)}Y^{\le -2}(u, y)}
\left(\frac{x_2+x_0}{x_1}\right)^{Y_0(u)}
\left(\frac{x_2+x_0}{x_1}\right)^{{\mathcal{L}}_{g}}v_1\nn
&= x_1^{-1}\delta\left(\frac{x_2+x_0}{x_1}\right)
x_{1}^{-Y_0(u)} e^{-\int_{0}^{-x_{1}}Y^{\le -2}(u, y)}
\left(\frac{x_2+x_0}{x_1}\right)^{{\mathcal{L}}_{g}}v_1\nn
&= x_1^{-1}\delta\left(\frac{x_2+x_0}{x_1}\right)
\left(\frac{x_2+x_0}{x_1}\right)^{{\mathcal{L}}_{g}}x_{1}
^{-Y_0(u)} e^{-\int_{0}^{-x_{1}}Y^{\le -2}(u, y)}v_1\nn
& =x_1^{-1}\delta\left(\frac{x_2+x_0}{x_1}\right)
\left(\frac{x_2+x_0}{x_1}\right)^{{\mathcal{L}}_g}
\Delta^{(u)}_V(x_1)v_1.
\end{align*}
\epfv

\begin{rema}\label{correction}
{\rm Here the authors would like to correct a sign mistake 
in \cite{H-log-twisted-mod}.
In \cite{H-log-twisted-mod}, $\Delta_{V}^{(u)}(x)$ is defined to be 
$x^{Y_{0}(u)}e^{\int_{0}{x}Y^{\le 2}(u, y)}$. But the correct definition is
$$\Delta_{V}^{(u)}(x)=x^{-Y_{0}(u)}e^{-\int_{0}^{-x}Y^{\le -2}(u, y)}.$$
In fact, using the definition with the wrong signs, in the proof of Theorem 5.5
in \cite{H-log-twisted-mod}, 
the author proved 
\begin{equation}\label{inverse-defn}
(Y_{W}^{(u)})^{e^{2\pi iY_{0}(u)}; p} (e^{2\pi iY_{0}(u)}v, z)w
= (Y_{W}^{(u)})^{e^{2\pi iY_{0}(u)}; p+1}(v, z)w.
\end{equation}
But equivarience property is formulated as
$$(Y_{W}^{(u)})^{e^{2\pi iY_{0}(u)}; p+1} (e^{2\pi iY_{0}(u)}v, z)w
= (Y_{W}^{(u)})^{e^{2\pi iY_{0}(u)}; p}(v, z)w$$
in the 
definition of twisted modules (see Condition 2 in Definition 3.1 
in \cite{H-log-twisted-mod}).
After the sign mistake is corrected, we obtain the  equivariance property in the definition. 
See the proof of Theorem \ref{main-thm}
above. Certainly, if we had used (\ref{inverse-defn}) as the equivariance 
property in the definition of generalized twisted module, then
$\Delta_{V}^{(u)}(x)$ given in \cite{H-log-twisted-mod} would be correct. }
\end{rema}

With additional conditions on $W$, we have the following consequence 
of Thereom \ref{main-thm}:

\begin{thm}\label{main-thm-1}
Let $V$, $g$ and $u$ be the same as in Theorem \ref{main-thm}. 
Let $(W, Y^{g}_{W})$
be a $\C/\Z$-graded 
generalized $g$-twisted $V$-module. 
Assume that $W$ is a direct sum  of generalized eigenspaces of 
$(Y_{W})_{0}(u)$.  
Then $(W,Y_W^{g g_u})$ 
is a $\C/\Z$-graded 
generalized $gg_u$-twisted $V$-module.
\end{thm}
\pf
First, we note that, by the commutator formula, we have 
$Y_1(u)\omega= u$ (we note here that \cite{H-log-twisted-mod} 
has a sign mistake in the proof of this fact and states $Y_1(u)\omega 
= -u$, but it should be $Y_1(u)\omega = u$ where $u \in V_{(1)}$ 
and $L(1)u = 0$). Using a proof identical to Lemma 5.3 in 
\cite{H-log-twisted-mod}, we have that 
the vertex operator $Y_{W}^{gg_u}(\omega, x)$ is in 
${\rm End}\;W[[x, x^{-1}]]$
and if we write 
$$Y_{W}^{gg_u}(\omega, x)=\sum_{n\in \Z}L_{W}^{gg_u}(n)x^{-n-2},$$
then 
\begin{equation}\label{new-L(0)}
L_{W}^{gg_u}(0)=L^{g}_{W}(0)
-(Y^g_W)_0(u)+\frac{1}{2} \mu
\end{equation}
where $L_{W}^{g}(0)$ is the corresponding Virasoro operator for the 
generalized $g$-twisted $V$-module $W$ and $\mu \in \C$ is
given by $Y_1(u)u = \mu {\bf 1}$. 
The eigenvalues of $e^{2 \pi i (Y^g_W)_0(u)}$ are of the form $e^{2 \pi i \beta}$ where 
$\beta \in P^{g_u}_W$. Also, a generalized eigenvector of $e^{2 \pi i (Y^g_W)_0(u)}$ 
with eigenvalue $e^{2 \pi i \beta}$ is also a generalized eigenvector of $(Y^g_W)_0(u)$ 
with eigenvalue $\beta + p$ for some $p \in \mathbb{Z}$. We note here that, by the 
commutator fomula, $[L^g_W(0),(Y^g_W)_0(u)] = 0$. By assumption, we have 
$$W = \coprod_{n \in \mathbb{C}, \alpha \in P^g_W} W^{[\alpha]}_{[n]},$$ 
where $W^{[\alpha]}_{[n]}$ is the homogeneous subspace of $W$ of generalized 
eigenvectors of $g$ with eigenvalue $e^{2 \pi i \alpha}$ and of $L^g_W(0)$ of 
eigenvalue $n$.  By assumption, we also have
$$W^{[\alpha]}_{[n]} =\coprod_{\tilde{\beta} \in P^{g_u}_W + \Z } 
W^{[\alpha],[\tilde{\beta}]}_{[n]}$$ 
for each $\alpha \in P^g_W$ and $n \in \C$. In fact all three gradings
are compatible 
since all the relevant operators commute. 
Let
$$ W_{\langle n \rangle}^{[\alpha],[\beta]} = \coprod_{ \tilde{\beta} 
\in \beta + \Z} W_{[n+\tilde{\beta}
-\frac{1}{2}\mu]}^{[\alpha],[\tilde{\beta}]}$$
for $n \in \C, \alpha \in P^g_W, \beta \in P^{g_u}_W$. The elements of 
$W^{[\alpha],[\beta]}_{\langle n \rangle}$ are generalized eigenvectors for
$L^{gg_u}_W(0)$ with eigenvalue $n$ (by  (\ref{new-L(0)})), 
generalized eigenvectors for 
$g$ with eigenvalue $e^{2 \pi i \alpha}$ and  generalized eigenvectors of
$g_u$ with eigenvalue $e^{2 \pi i \beta}$. In particular, they are 
generalized eigenvectors of 
$gg_{u}$ with eigenvalue $e^{2 \pi i (\alpha+\beta)}$.
We have
$$ W = \coprod_{n \in \C, \alpha \in P^g_W, \beta \in P^{g_u}_W } 
W^{[\alpha],[\beta]}_{\langle n \rangle}.$$
Since all operators commute, we have 
$$ P^{gg_u}_W = \{ \gamma \in \mathbb{C} | \Re(\gamma) \in [0,1), \gamma 
\in \alpha + \beta + \Z, \text{ for } \alpha \in P^g_W, \beta \in P^{g_u}_W\}.$$
Define
$$ W^{\langle\gamma\rangle}_{\langle n \rangle} 
= \coprod_{\substack{\alpha \in P^g_W, 
\beta \in P^{g_u}_W \\ \alpha + \beta \in \gamma +\Z}}
W^{[\alpha],[\beta]}_{\langle n \rangle}$$
Then $W$ is doubly graded by the eigenvalues of 
$L^{gg_u}_W(0)$ and $gg_u$, 
that is, $W^{\langle \gamma \rangle}_{\langle n \rangle}$ 
is the intersection of the generalized 
eigenspace of $L^{gg_u}(0)$ with eigenvalue $n$ and the 
generalized eigenspace of $gg_u$ 
with eigenvalue $e^{2 \pi i \gamma}$, and 
$$ W = \coprod_{n \in \C, \gamma \in P^{gg_h}_W} 
W^{\langle \gamma \rangle}_{\langle n \rangle}.$$
Together with Theorem \ref{main-thm}, we see that 
 $W$ is a $\C/\Z$-graded generalized $gg_u$-twisted module.
\epfv

\begin{rema}
{\rm In this paper, $V$ is always a vertex operator algebra. In particular, 
$V$ has a conformal vector $\omega$ such that the gradings on $V$ and on 
generalized twisted $V$-modules 
are given by the coefficients of $x^{-2}$ of the vertex operators 
of the conformal element $\omega$. 
But Theorem \ref{main-thm} also holds in the case that $V$ is a 
grading-restricted vertex algebra (that is,
a vertex algebra with grading-restricted grading given by $L(0)$) 
since we do not use anything involving 
the conformal element. 
Theorem \ref{main-thm-1} can be generalized to the case that 
$V$ is a grading-restricted vertex algebra
(likely satisfying some additional conditions) by defining 
\begin{align*}
Y^{gg_u}_W(v,x) &= Y^g_W(\Delta^{(u)}_V(x) v, x),\\
L^{gg_u}_{W}(0)&=L^{g}_{W}(0)
-(Y^g_W)_0(u)+\frac{1}{2} \mu,\\
L^{g g_u}_{W}(-1)&=L^{g}_{W}(-1)-(Y^{g}_{W})_{-1}(u)
\end{align*}
for $v\in V$, where $\mu\in \C$ is given by $(Y_{V})_{1}(u)u=\mu \one$.}
\end{rema}

The generalized $gg_u$-twisted $V$-module in Theorem \ref{main-thm} is $\C/\Z$-graded. 
We now discuss $\C$-graded generalized $gg_u$-twisted modules. 
Note that 
$g=e^{2\pi i (\mathcal{S}_{g}+\mathcal{N}_{g})}$ and $g_{u}=e^{2\pi i Y_{0}(u)}$. 
If there is another semisimple operator $\widetilde{\mathcal{S}}_{g}$ on $V$ such that 
$g=e^{2\pi i (\widetilde{\mathcal{S}}_{g}+\mathcal{N}_{g})}$, then 
the set of eigenvalues of $\widetilde{\mathcal{S}}_{g}$ must be contained in 
$P_{V}^{g}+\Z$. In particular, we have
$$V = \coprod_{\tilde{\alpha}\in P_{V}^{g}+\Z}V^{[\tilde{\alpha}]}
=\coprod_{\tilde{\alpha}\in P_{V}^{g}+\Z,
n\in \Z}V_{(n)}^{[\tilde{\alpha}]},$$
where  $V_{(n)}^{[\tilde{\alpha}]}$
for $\tilde{\alpha}\in P_{V}^{g}+\Z$  and $n\in \Z$ is the eigenspace 
of $\widetilde{\mathcal{S}}_{g}$ restricted $V_{(n)}$ with eigenvalue $\tilde{\alpha}$
 and 
$$V^{[\tilde{\alpha}]}=\coprod_{n\in \Z}V_{(n)}^{[\tilde{\alpha}]}$$ 
for $\tilde{\alpha}\in P_{V}^{g}+\Z$ is the eigenspace 
of $\widetilde{\mathcal{S}}_{g}$ with eigenvalue $\tilde{\alpha}$.

We now assume that there is a semisimple operator $\widetilde{\mathcal{S}}_{g}$ on $V$
such that $g=e^{2\pi i (\widetilde{\mathcal{S}}_{g}+\mathcal{N}_{g})}$ and 
$Y_{V}(u, x)v\in V^{[\tilde{\alpha}_{1}+\tilde{\alpha}_{2}]}$ for $u\in V^{[\tilde{\alpha}_{1}]}$ 
and $v\in V^{[\tilde{\alpha}_{2}]}$. 
In this case, 
$$gg_{u}=e^{2\pi i (\widetilde{\mathcal{S}}_{g}+\mathcal{N}_{g})}e^{2\pi i Y_{0}(u)}
=e^{2\pi i (\widetilde{\mathcal{S}}_{g}+Y_{0}(u)_{S})}
e^{2\pi i\mathcal{N}_{g}}e^{2\pi i (Y_{0}(u)_N)}.$$
In particular, 
$e^{2\pi i (\widetilde{\mathcal{S}}_{g}+Y_{0}(u)_{S})}$ is the semisimple part of 
$gg_{u}$. We denote the operator 
$\widetilde{\mathcal{S}}_{g}+Y_{0}(u)_{S}$ by 
$\widetilde{\mathcal{S}}_{gg_{u}}$. 

Let $W$ be a generalized $g$-twisted $V$-module with a semisimple 
action of $\widetilde{\mathcal{S}}_{g}$
such that the actions of $e^{2\pi i\mathcal{S}_{g}}$  and 
$e^{2\pi i\widetilde{\mathcal{S}}_{g}}$
on $W$ are equal.
Then $W$ is a direct sum of eigenspaces of $\widetilde{\mathcal{S}}_{g}$
and the set of eigenvalues of 
$\widetilde{\mathcal{S}}_{g}$ on $W$ must be contained in 
$P_{W}^{g}+\Z$. In particular, we have a grading
$$W = \coprod_{\widetilde{\alpha}\in P_{W}^{g}+\Z}W^{[\tilde{\alpha}]}
=\coprod_{\tilde{\alpha}\in P_{W}^{g}+\Z,
n\in \C}W_{[n]}^{[\tilde{\alpha}]},$$
where $W_{[n]}^{[\tilde{\alpha}]}$
for $\tilde{\alpha}\in P_{V}^{g}+\Z$  and $n\in \C$ is the eigenspace 
of $\widetilde{\mathcal{S}}_{g}$ restricted $W_{[n]}$ with eigenvalue $\tilde{\alpha}$ and
$$W^{[\tilde{\alpha}]}=\coprod_{n\in \C}W_{[n]}^{[\tilde{\alpha}]}$$ 
for $\tilde{\alpha}\in P_{W}^{g}+\Z$
is the eigenspace
of $\widetilde{\mathcal{S}}_{g}$  with eigenvalue $\tilde{\alpha}$.
We further assume that this grading is compatible with the twisted modue structure, that is, 
$Y^{g}_{W}(v, x)w\in W^{[\tilde{\alpha}_{1}+\tilde{\alpha}_{2}]}$ for $u\in 
V^{[\tilde{\alpha}_{1}]}$ and $w\in W^{[\tilde{\alpha}_{2}]}$. 

Since $V$ is a vertex operator algebra, $V^{[\tilde{\alpha}]}_{(n)}$ 
for $\tilde{\alpha}\in  P_{W}^{g}+\Z$ and $n\in \Z$
is finite dimensional. 
In particular, it can be decomposed into generalized 
$Y_0(u)$-eigenspace
$$V_{(n)}^{[\tilde{\alpha}]}  = \coprod_{\tilde{\beta} \in P_{V}^{g_{u}}+\Z} 
V_{(n)}^{[\tilde{\alpha}],[\tilde{\beta}]},$$
where $V_{(n)}^{[\tilde{\alpha}],[\tilde{\beta}]}$ is the generalized eigenspace of $Y_0(u)$ restricted to 
$V_{(n)}^{[\tilde{\alpha}]}$ with eigenvalue $\tilde{\beta}$.
By assumption, $W^{[\tilde{\alpha}]}_{[n]}$ can be decomposed into generalized 
$(Y^g_W)_0(u)$-eigenspaces, that is, 
$$W_{[n]}^{[\tilde{\alpha}]}  = \coprod_{\tilde{\beta} \in P_{W}^{g_{u}}+\Z} 
W_{[n]}^{[\tilde{\alpha}],[\tilde{\beta}]},$$
where $W_{[n]}^{[\tilde{\alpha}],[\tilde{\beta}]}$ is the generalized eigenspace of 
$(Y^g_W)_0(u)$ restricted 
to $W_{[n]}^{[\tilde{\alpha}]}$ with eigenvalue $\tilde{\beta}$. 

On $V_{(n)}^{[\tilde{\alpha}],[\tilde{\beta}]}$, $\widetilde{\mathcal{S}}_{gg_{u}}$ acts as 
$\tilde{\alpha}+\tilde{\beta}$ so that $V_{(n)}^{[\tilde{\alpha}],[\tilde{\beta}]}$ is in 
the generalized eigenspace of  $\widetilde{\mathcal{S}}_{gg_{u}}$ with eigenvalue 
$\tilde{\alpha}+\tilde{\beta}$. For 
$\tilde{\gamma}\in P_{V}^{gg_{u}}+\Z$ and $n\in \Z$, 
let 
\begin{align*}
V^{\langle \tilde{\gamma} \rangle}_{(n)} 
= \coprod_{\substack{\tilde{\alpha}\in P_{V}^{g}+\Z, \tilde{\beta} \in P_{V}^{g_{u}}+\Z\\ 
\tilde{\alpha}+\tilde{\beta} = \tilde{\gamma}}} V_{(n)}^{[\tilde{\alpha}],[\tilde{\beta}]}
.
\end{align*}
Then $V^{\langle \tilde{\gamma} \rangle}_{(n)}$ is the generalized eigenspace of 
$gg_{h}$ restricted to $V_{(n)}$ with eigenvalue $\tilde{\gamma}$.
Similarly,  on $W_{[n]}^{[\tilde{\alpha}],[\tilde{\beta}]}$, $\widetilde{\mathcal{S}}_{gg_{u}}$ acts as 
$\tilde{\alpha}+\tilde{\beta}$ so that $W_{[n]}^{[\tilde{\alpha}],[\tilde{\beta}]}$ is in 
the generalized eigenspace of  $\widetilde{\mathcal{S}}_{gg_{u}}$ with eigenvalue 
$\tilde{\alpha}+\tilde{\beta}$. For $\tilde{\gamma}\in P_{W}^{gg_{u}}+\Z$ and $n\in \C$, 
let 
\begin{align*}
W^{\langle \tilde{\gamma} \rangle}_{\langle n\rangle} 
= \coprod_{\substack{\tilde{\alpha}\in P_{W}^{g}+\Z, \tilde{\beta} 
\in P_{W}^{g_{u}}+\Z\\ 
\tilde{\alpha}+\tilde{\beta} = \tilde{\gamma}}} W_{[n+\tilde{\beta}
-\frac{1}{2}\mu]}^{[\tilde{\alpha}],[\tilde{\beta}]}
.
\end{align*}
Then $W^{\langle \tilde{\gamma} \rangle}_{\langle n\rangle}$ is the 
generalized eigenspace of 
$gg_{u}$ restricted to 
\begin{equation}\label{ggh-grading}
W_{\langle n\rangle}=\coprod_{\tilde{\gamma} \in P_{W}^{gg_{u}}+\Z}
W^{\langle \tilde{\gamma} \rangle}_{\langle n\rangle}
\end{equation}
with eigenvalue $\tilde{\gamma}$.
Then these subspaces give new gradings to $V$ and $W$, that is,
$$V=\coprod_{n\in \Z, \tilde{\gamma}\in P_{V}^{gg_{u}}+\Z}
V_{(n)}^{\langle \tilde{\gamma}\rangle}$$
and 
$$W=\coprod_{n\in \C, \tilde{\gamma}\in P_{W}^{gg_{u}}+\Z}
W_{\langle n\rangle}^{\langle \tilde{\gamma}\rangle}.$$

\begin{thm}\label{C-graded}
Let $V$, $g$ and $u$ be the same as in Theorem \ref{main-thm}. 
Let $W$ be a generalized $g$-twisted $V$-module. Assuming that 
there is a semisimple operator $\widetilde{\mathcal{S}}_{g}$ on $V$ 
such that $g=e^{2\pi i (\widetilde{\mathcal{S}}_{g}+\mathcal{N}_{g})}$ and 
$Y_{V}(u, x)v\in V^{[\tilde{\alpha}_{1}+\tilde{\alpha}_{2}]}$ for $u\in V^{[\tilde{\alpha}_{1}]}$ 
and $v\in V^{[\tilde{\alpha}_{2}]}$, where for $\alpha\in P_{V}^{g}$, $V^{[\tilde{\alpha}]}$ is the 
eigenspace of $\widetilde{\mathcal{S}}_{g}$ with eigenvalue $\tilde{\alpha}$. 
Assume also that $\widetilde{\mathcal{S}}_{g}$ acts on $W$ semisimply, 
 the actions of $e^{2\pi i\mathcal{S}_{g}}$  and $e^{2\pi i\widetilde{\mathcal{S}}_{g}}$
on $W$ are equal
and 
$Y^{g}_{W}(v, x)w\in W^{[\tilde{\alpha}_{1}+\tilde{\alpha}_{2}]}$ for $u\in 
V^{[\tilde{\alpha}_{1}]}$ and $w\in W^{[\tilde{\alpha}_{2}]}$.
Then $(W,Y^{gg_u}_W)$ with the gradings of $W$ given by 
(\ref{ggh-grading})
 is a $\C$-graded generalized $gg_u$-twisted module.
In particular, if we  define a $\C/\Z$-grading of $W$  by 
$$W=\coprod_{\alpha\in P_{W}^{g}}W^{(\alpha)},$$
where for $\alpha\in P_{W}^{g}$, 
$$W^{(\alpha)}=\coprod_{k\in \Z}W^{[\alpha+k]},$$
then $(W,Y^{gg_u}_W)$ with this $\C/\Z$-grading is a $\C/\Z$-graded generalized $gg_u$-twisted module.
\end{thm}
\pf
We need only show that the new gradings on $V$ and $W$ are compatible with the 
twisted vertex operator map. Indeed, suppose 
$v \in V_{(n)}^{[\tilde{\alpha}],[\tilde{\beta}]}\subset 
V^{\langle \tilde{\gamma} \rangle}_{(n)}$, where  
$\tilde{\alpha}\in P_{V}^{g}+\Z$, $\tilde{\beta} \in P_{V}^{g_{u}}+\Z$ and 
$\tilde{\alpha}+\tilde{\beta} = \tilde{\gamma}$, 
and $w \in W_{[n+\tilde{\beta}'-\frac{1}{2}\mu]}^{[\tilde{\alpha}'],[\tilde{\beta}']}\subset W^{\langle \tilde{\gamma} \rangle}_{\langle n\rangle} $,
where $\tilde{\alpha}'\in P_{W}^{g}+\Z$, $\tilde{\beta}' \in P_{W}^{g_{u}}+\Z$ and
$\tilde{\alpha}'+\tilde{\beta}' = \tilde{\gamma}'$. 
Then 
\begin{align*}
&\widetilde{\mathcal{S}}_{gg_{h}}Y_{W}^{gg_{u}}(v,x)w\nn
&\quad =(\widetilde{\mathcal{S}}_{g}+Y_{0}(u)_{S})Y_{W}^{g}(\Delta_{V}^{(u)}(x)v,x)w \nn
&\quad=\widetilde{\mathcal{S}}_{g}Y_{W}^{g}(\Delta_{V}^{(u)}(x)v,x)w+Y_{0}(u)_{S}
Y_{W}^{g}(\Delta_{V}^{(u)}(x)v,x)w \nn
&\quad=(\tilde{\alpha}+\tilde{\alpha}')Y_{W}^{g}(\Delta_{V}^{(u)}(x)v,x)w+
Y_{W}^{g}(Y_{0}(u)_{S}\Delta_{V}^{(u)}(x)v,x)w +Y_{W}^{g}(\Delta_{V}^{(u)}(x)v,x)Y_{0}(u)_{S}w\nn
&\quad=(\tilde{\alpha}+\tilde{\alpha}')Y_{W}^{g}(\Delta_{V}^{(u)}(x)v,x)w+
Y_{W}^{g}(\Delta_{V}^{(u)}(x)Y_{0}(u)_{S}v,x)w +Y_{W}^{g}(\Delta_{V}^{(u)}(x)v,x)\tilde{\beta}'w\nn
&\quad=(\tilde{\alpha}+\tilde{\alpha}'+\tilde{\beta}+\tilde{\beta}')Y_{W}^{gg_{u}}(v,x)w\nn
&\quad=(\tilde{\gamma}+\tilde{\gamma}')Y_{W}^{gg_{u}}(v,x)w.
\end{align*}
This means 
$$Y_{W}^{gg_{u}}(v,x)w
\in W^{\langle \tilde{\gamma}+\tilde{\gamma}'\rangle} \{ x \} [\log x],$$
proving that the new gradings on $V$ and $W$ are indeed compatible with the 
twisted vertex operator map.
\epfv

We now discuss the case of strongly $\C$-graded generalized $g$-twisted $V$-modules. 
Let $W$ be a $\C$-graded generalized $g$-twisted $V$-module such that the assumptions in 
Theorem \ref{C-graded} hold. Then $W$ with the 
$\C/\Z$-grading
$$W=\coprod_{\alpha\in P_{W}^{g}}W^{(\alpha)}=\coprod_{\alpha\in P_{W}^{g},
k\in \Z}W^{[\alpha+k]}$$
is a $\C/\Z$-graded generalized $g$-twisted $V$module. 
For each $\alpha \in P^g_W$ and $n\in \C$,
we denote the set consisting of those $\tilde{\alpha}\in \alpha+\Z$ such that 
$W_{[n]}^{[\tilde{\alpha}]}$ is not $0$ by $Q^{\alpha}_{n}$. 
Let $Q_{n}=\bigcup_{\alpha\in P_{W}^{g}}Q_{n}^{\alpha}$ for $n\in \C$ and $Q=\bigcup_{n\in \C}Q_{n}$.

If $W$ is strongly $\C/\Z$-graded, then
$$W_{[n]}^{(\alpha)}=\coprod_{\alpha\in P_{W}^{g},
\tilde{\alpha}\in \alpha+\Z}W_{[n]}^{[\tilde{\alpha}]}$$ 
for $\alpha\in P_{W}^{g}$ and $n\in \C$ are finite dimensional and for each $\alpha\in P_{W}^{g}$ and $n\in \C$,
there are only finitely many 
$\tilde{\alpha}\in \alpha+\Z$ such that $W_{[n]}^{[\tilde{\alpha}]}$ is not $0$.
In particular, for each $\alpha\in P_{W}^{g}$ and $n\in \C$, $Q^{\alpha}_{n}$ is a finite set
in this case.

\begin{thm}\label{strongly-C-graded}
Let $V$, $g$ and $u$ be the same as in Theorem \ref{main-thm}. 
Let $W$ be a $\C$-graded generalized $g$-twisted $V$-module such that the assumptions in 
Theorem \ref{C-graded} hold. Assume in addtion 
that the corresponding $\C/\Z$-graded 
 generalized $g$-twisted $V$-module is strongly $\C/\Z$-graded, $Q$ 
is a finite set, 
then the $\mathbb{C}$-graded generalized $gg_u$-twisted module 
$(W,Y_W^{gg_u})$ is  strongly $\mathbb{C}$-graded.
\end{thm}
\pf
For $\tilde{\gamma}\in P_{W}^{gg_{u}}$ and $n\in \C$, 
\begin{align}\label{direct-sum}
W^{\langle \tilde{\gamma} \rangle}_{\langle n\rangle} 
= \coprod_{\substack{\tilde{\alpha}\in P_{W}^{g}+\Z, \tilde{\beta} \in P_{W}^{g_{u}}+\Z\\ 
\tilde{\alpha}+\tilde{\beta} = \tilde{\gamma}}} 
W_{[n+\tilde{\beta}-\frac{1}{2}\mu]}^{[\tilde{\alpha}],[\tilde{\beta}]}
=\coprod_{\substack{\tilde{\alpha}\in Q_{n+\tilde{\beta}-\frac{1}{2}\mu}, \;
\tilde{\beta} \in P_{W}^{g_{u}}+\Z\\ 
\tilde{\alpha}+\tilde{\beta} = \tilde{\gamma}}} 
W_{[n+\tilde{\beta}-\frac{1}{2}\mu]}^{[\tilde{\alpha}],[\tilde{\beta}]}
=\coprod_{\substack{\tilde{\alpha}\in Q, \;
\tilde{\beta} \in P_{W}^{g_{u}}+\Z\\ 
\tilde{\alpha}+\tilde{\beta} = \tilde{\gamma}}} 
W_{[n+\tilde{\beta}-\frac{1}{2}\mu]}^{[\tilde{\alpha}],[\tilde{\beta}]}.
\end{align}
Since $Q$ is a finite set,  the right hand side of (\ref{direct-sum}) is  a finite direct sum. 
Since in this case, $W_{[n+\tilde{\beta}-\frac{1}{2}\mu]}^{[\tilde{\alpha}],[\tilde{\beta}]}$
are all finite dimensional, 
we see that $W^{\langle \tilde{\gamma} \rangle}_{\langle n\rangle}$ are all finite dimensional. 
 
Since $W$  is strongly $\C/\Z$-graded, for each 
$n\in \C$, $\tilde{\alpha}\in P_{W}^{g}+\Z$, $\tilde{\beta}\in P_{W}^{g_{u}}$ and 
$\tilde{\gamma}\in P_{W}^{gg_{u}}$, 
there exists $L_{n, \tilde{\alpha}, \tilde{\beta}, \tilde{\gamma}}\in \Z$ such that 
$W_{[n-l+\tilde{\beta}-\frac{1}{2}\mu]}^{[\tilde{\alpha}],[\tilde{\beta}]}=0$
when $l\in \Z$ is larger than $L_{n, \tilde{\alpha}, \tilde{\beta}, \tilde{\gamma}}$.
Let $L$ be an integer larger than $L_{n, \tilde{\alpha}, \tilde{\beta}, \tilde{\gamma}}$
for $n\in \C$, $\tilde{\alpha}\in P_{W}^{g}+\Z$, $\tilde{\beta}\in P_{W}^{g_{u}}$ and 
$\tilde{\gamma}\in P_{W}^{gg_{u}}$
satisfying $\tilde{\alpha}+\tilde{\beta} = \tilde{\gamma}$. Then 
$$W^{\langle \tilde{\gamma} \rangle}_{\langle n-l\rangle} 
=\coprod_{\substack{\tilde{\alpha}\in Q, \;
\tilde{\beta} \in P_{W}^{g_{u}}+\Z\\ 
\tilde{\alpha}+\tilde{\beta} = \tilde{\gamma}}} 
W_{[n-l+\tilde{\beta}-\frac{1}{2}\mu]}^{[\tilde{\alpha}],[\tilde{\beta}]}=0.$$
\epfv

Since $\Delta^{(u)}_V(x)$ is invertible, we immediately have:

\begin{cor}
Let $V$, $g$ and $u$ be the same as in Theorem \ref{main-thm}. 
If $(W,Y^g_W)$ is an irreducible generalized $g$-twisted $V$-module, 
then $(W,Y^{gg_u}_W)$ is an irreducible generalized $gg_u$-twisted $V$-module.
\end{cor}

These results in fact give us functors from suitable subcategories of 
generalized twisted $V$-modules to themselves. 
Let $\mathcal{C}$ be the category of all generalized twisted
$V$-modules. Let $u\in V_{(1)}$ and $G^{u}$ the subgroup of the 
automorphism group of $V$ consisting the automorphisms
$g$ fixing $u$ (that is, $g(u)=u$). Consider the subcategory $\mathcal{C}^{u}$ 
of $\mathcal{C}$ consisting
of generalized $g$-twisted $V$-mdules for $g\in G^{u}$.
For $g\in G^{u}$, we also have $gg_{u}\in  G^{u}$.
Then we have a functor 
$$\Delta^{u}: \mathcal{C}^{u}\to \mathcal{C}^{u}$$
defined as follows: For an object 
$(W, Y_{W}^{g})$ of 
$\mathcal{C}^{u}$ (a generalized $g$-twisted $V$-module for $g\in G^{u}$),  
$$\Delta^{u}(W, Y_{W}^{g})=(W, Y^{gg_u}_W).$$
For a morphism (a module map) $f$
from an object $(W_{1}, Y_{W_{1}}^{g})$
of $\mathcal{C}^{u}$ to another 
object $(W_{2}, Y_{W_{2}}^{g})$ of 
$\mathcal{C}^{u}$, we have 
\begin{align*}
f(Y_{W_{1}}^{gg_{u}}(v,x)w_{1})&
=f(Y_{W_{1}}^{g}(\Delta_{V}^{(u)}(x)v,x)w_{1})\nn
&=Y_{W_{2}}^{g}(\Delta_{V}^{(u)}(x)v,x)f(w_{1})\nn
&=Y_{W_{2}}^{gg_{u}}(v,x)f(w_{1}).
\end{align*}
Hence $f$ is also a module map from 
$(W_{1}, Y_{W_{1}}^{gg_{u}})$ to 
$(W_{2}, Y_{W_{2}}^{gg_{u}})$. We define 
the image of $f$ under $\Delta^{u}$ to be $f$ viewed as a 
module map from 
$(W_{1}, Y_{W_{1}}^{gg_{u}})$ to 
$(W_{2}, Y_{W_{2}}^{gg_{u}})$. 

We have the following result:

\begin{thm}
For $u\in V_{(1)}$, $\Delta^{u}$ is a functor from 
$\mathcal{C}^{u}$ to itself. Moreover, $\Delta^{u}$ 
is in fact an automorphism of the category $\mathcal{C}^{u}$. 
\end{thm}
\pf
It is clear that $\Delta^{u}$ is a functor. 
Note that $\mathcal{C}^{-u}=\mathcal{C}^{u}$.
Then we also have a functor $\Delta^{-u}$
 from $\mathcal{C}^{u}$ to itself. By Proposition \ref{Delta-prop},
 we see that 
 $$\Delta^{u}\circ \Delta^{-u}=\Delta^{-u}\circ \Delta^{u}
 =1_{u},$$
 where $1_u$ is the identity functor on 
 $\mathcal{C}^{u}$. 
 \epfv
 
 \begin{rema}
 {\rm Since the underlying vector space of 
 $(W, Y_{W}^{g})$ and $\Delta^{u}(W, Y_{W}^{g})
 =(W, Y^{gg_u}_W)$ are the same, 
 $\Delta^{u}$ is indeed an automorphism of the 
 category $\mathcal{C}^{u}$ (that is, an isomorphism from 
$\mathcal{C}^{u}$ to itself), not just an equivalence from 
 $\mathcal{C}^{u}$ to itself.}
 \end{rema}

\section{Affine vertex (operator) algebras and their automorphisms}

In the next two sections, we shall apply the automorphisms 
of categories obtained in the preceding section to construct certain particular
generalized twisted modules for affine vertex operator
algebras. In this section,
we recall some basic facts on affine vertex operator algebras
and their automorphisms. 

Let $\g$ be a finite-dimensional simple Lie algebra and 
$(\cdot, \cdot)$ its normalized 
Killing form. 
Let $g$ be an automorphism of $\g$. Assume also that 
$(\cdot, \cdot)$ is invariant under $g$. 
Then $g=e^{2\pi i\mathcal{S}_{g}} e^{2\pi i \mathcal{N}_{g}}$,
where $e^{2\pi i\mathcal{S}_{g}}$ and $ e^{2\pi i \mathcal{N}_{g}}$ 
are the semisimple and unipotent parts
of $g$. Let 
$$P_{\g}^g=\{\alpha\in \C\;|\; \Re(\alpha)\in [0, 1), 
e^{2 \pi i\alpha} \;\text{is an eigenvalue of}\;g\}.$$
Then 
$$\mathfrak{g}=\coprod_{\alpha\in P_{\g}}
\mathfrak{g}^{[\alpha]},$$
where for $\alpha\in P_{\g}$,
$\mathfrak{g}^{[\alpha]}$ is the generalized eigenspace of $g$ 
(or the eigenspace of 
$e^{2\pi i\mathcal{S}_{g}}$) with the eigenvalue $e^{2\pi i\alpha}$.

By Proposition 5.3 in \cite{H-aff-va-twisted-mod},
$e^{2\pi i\mathcal{S}_{g}}=\sigma=\tau_{\sigma}\mu 
e^{2\pi i {\rm ad}_{h}}\tau_{\sigma}^{-1}$ where $h \in \mathfrak{h}$, 
$\mu$ is a diagram automorphism of $\mathfrak{g}$ such that
$\mu(h) = h$ and $\tau_\sigma$ is an automorphism of $\mathfrak{g}$. 
By Proposition 5.4 in \cite{H-aff-va-twisted-mod}, $\mathcal{N}_{g}
={\rm ad}_{a_{\mathcal{N}_{g}}}$, where 
$a_{\mathcal{N}_{g}}\in \g^{[0]}$. Then 
$$\tau_{\sigma}^{-1}\mathcal{N}_{g}\tau_{\sigma}
=\tau_{\sigma}^{-1}{\rm ad}_{a_{\mathcal{N}_{g}}}\tau_{\sigma}
={\rm ad}_{\tau_{\sigma}^{-1}a_{\mathcal{N}_{g}}}.$$
So 
\begin{align*}
g&=e^{2\pi i\mathcal{S}_{g}}e^{2\pi i\mathcal{N}_{g}}
=\tau_{\sigma}\mu e^{2\pi i {\rm ad}_{h}}\tau_{\sigma}^{-1}e^{2\pi i\mathcal{N}_{g}}\nn
&=\tau_{\sigma}\mu e^{2\pi i {\rm ad}_{h}}e^{2\pi i\tau_{\sigma}^{-1}\mathcal{N}_{g}\tau_{\sigma}}
\tau_{\sigma}^{-1}
=\tau_{\sigma}\mu e^{2\pi i {\rm ad}_{h}}e^{2\pi i{\rm ad}_{\tau_{\sigma}^{-1}a_{\mathcal{N}_{g}}}}
\tau_{\sigma}^{-1}.
\end{align*}
Let $g_{\sigma}=\mu e^{2\pi i {\rm ad}_{h}}e^{2\pi i{\rm ad}_{\tau_{\sigma}^{-1}a_{\mathcal{N}_{g}}}}$. Then
$\mu e^{2\pi i {\rm ad}_{h}}$ and $e^{2\pi i{\rm ad}_{\tau_{\sigma}^{-1}a_{\mathcal{N}_{g}}}}$
are the semisimple and unipotent, respectively,  parts of $g_{\sigma}$. We shall denote 
the semisimple part $\mu e^{2\pi i {\rm ad}_{h}}$ of $g_{\sigma}$ by $g_{\sigma, s}$. 

Let $\widehat{\g} = \g \otimes \C[t,t^{-1}] \oplus \C {\bf k}$ equipped with the bracket operation
\begin{align*}
[a \otimes t^m, b \otimes t^n] &= [a,b] \otimes t^{m+n} + (a,b) m \delta_{m+n,0} {\bf k},\\
[a \otimes t^m, {\bf k}] &= 0
\end{align*}
for $a,b \in \g$ and $m,n \in \Z$. Let $\widehat{\g}_{\pm} = \g \otimes t^{\pm 1} \C[t^{\pm 1}]$ so that 
\begin{align*}
\widehat{\g} = \widehat{\g}_{-} \oplus \g \oplus \C {\bf k} \oplus \widehat{\g}_+.
\end{align*}
Consider $\C$ as a $1$-dimensional $\g \oplus \C{\bf k} \oplus \widehat{\g}_{+}$-module where $\g \oplus \widehat{\g}_+$ acts trivially, ${\bf k}$ acts as $\ell \in \C$,
and let
\begin{equation*}
M(\ell,0) = U(\widehat{\g}) \otimes_{U(\g \oplus \C{\bf k} \oplus \widehat{\g}_+)} \C
\end{equation*}
be the induced $\widehat{\g}$-module. Let $J(\ell,0)$ be the maximal proper $\widehat{\g}$-submodule of $M(\ell,0)$ and let $L(\ell,0) = M(\ell,0)/J(\ell,0)$.
It is well known (\cite{FZ}) that $M(\ell,0)$ and $L(\ell,0)$ have the structure of a vertex operator algebra when $\ell\ne -h^{\vee}$, where
$h^{\vee}$ is the dual Coxeter number of $\g$. 

The automorphisms $\mu$, $e^{2\pi i {\rm ad}_{h}}$, $e^{2\pi i{\rm ad}_{\tau_{\sigma}^{-1}a_{\mathcal{N}_{g}}}}$,
$g_{\sigma}$ and $g_{\sigma, s}$ also give automorphisms of $M(\ell, 0)$ and $L(\ell, 0)$. But on 
$M(\ell, 0)$ and $L(\ell, 0)$,
${\rm ad}_{h}$ and ${\rm ad}_{\tau_{\sigma}^{-1}a_{\mathcal{N}_{g}}}$ act as $h(0)$ and 
$(\tau_{\sigma}^{-1}a_{\mathcal{N}_{g}})(0)$. So the corresponding automorphisms of $L(\ell, 0)$
are actually $e^{2\pi i h(0)}$ and $e^{2\pi i (\tau_{\sigma}^{-1}a_{\mathcal{N}_{g}})(0)}$.

We will also study inner automorphisms of $\g$ below. Fix an element $a \in \g$. The element $a$ has a Chevalley-Jordan decomposition given by $a = s + n$ where ad$_s$ is semisimple on $\g$, ad$_n$ is nilpotent on $\g$, and $[s,n]=0$. In particular, if ad$_s$ is semisimple, then $s$ belongs to some Cartan subalgebra of $\g$. Without loss of generality, we take $s \in \mathfrak{h}$. Let $a_1,\dots, a_n \in \g$ be eigenvectors of ad$_s$ with eigenvalues $\lambda_1,\dots,\lambda_n$, respectively. That is, we have
$$
\mathrm{ad}_s(a_i) = [s,a_i]= \lambda_i a_i
$$
for $1 \le i \le n$. We note here that if $\lambda_i \neq 0$, we have that 
\begin{align}
(s,a_i) = \frac{1}{\lambda_i} (s,[s,a_i]) = \frac{1}{\lambda_i} ([s,s],a_i) = 0.
\end{align}
Hence, $(s,a_i) \neq 0$ is only possible when $\lambda_i=0$.

We now use $V$ to denote $M(\ell,0)$ or $L(\ell,0)$.
Then we have $Y_{V}(a(-1)\one,x) = \sum_{n \in \Z} a(n)x^{-n-1}$.
The action ad$_a$ on $\g$ becomes the action of $a(0)$ on $V$. In particular, we have that $a(0)$ is a derivation of $V$ since by the commutator formula for a vertex operator algebra we have
\begin{align}
Y(a(0)v,x) = [a(0),Y(v,x)]
\end{align}
for all $v \in V$. 
We note that, similarly, $s(0)$ and $n(0)$ are also derivations of $V$.
We have that $V$ is spanned by elements of the form
\begin{align}\label{general-eigenvector}
a_{i_1}(m_1) \cdots a_{i_k}(m_k) \one
\end{align}
where $k \ge 0$, $i_1,\dots,i_k \in \{1,\dots,n\}$ and $m_1,\dots, m_k \in \mathbb{Z}$. In general, elements of the form (\ref{general-eigenvector}) are eigenvectors of $s(0)$ and generalized eigenvectors of $a(0)$. That is, we have
\begin{align}
s(0) a_{i_1}(m_1) \cdots a_{i_k}(m_k) \one = (\lambda_{i_1} + \cdots + \lambda_{i_k}) a_{i_1}(m_1) \cdots a_{i_k}(m_k) \one,
\end{align}
We have that $g_a = e^{2 \pi i a(0)}$, an automorphism of $V$, and its semisimple and unipotent parts are $g_s = e^{2 \pi i s(0)}$ and $g_{\mathcal{U}} = e^{2 \pi i n(0)}$, respectively. We have that (\ref{general-eigenvector}) is an eigenvector for the automorphism $e^{2 \pi i s(0)}$ with eigenvalue $e^{2 \pi i (\lambda_{i_1} + \cdots + \lambda_{i_k})}$.  Thus, the actions of $a(0)$ and $s(0)$ give us a natural $\C$-grading on $V$, given by the eigenspaces of $s(0)$, which are also generalized eigenspaces of $a(0)$:
\begin{align}
V = \coprod_{\alpha \in P^g_V + \Z} V^{[\alpha]}
\end{align} 
where
\begin{align}
V^{[\alpha]} = \mathrm{span} \{ a_{i_1}(m_1)\cdots a_{i_k}(m_k) \one | k \ge 0, m_1, \dots, m_k \in \Z, \lambda_{i_1} + \cdots \lambda_{i_k} = \alpha\}.
\end{align}
Moreover, it is easy to see that the grading compatibility condition is satisfied by $V$: Let $v_1 \in V^{[\alpha]}$ and $v_2 \in V^{[\beta]}$ for some $\alpha, \beta \in \C$. Then, we have that
\begin{align}
s(0)Y(v_1,x)v_2 = Y(s(0)v_1,x)v_2 + Y(v_1,x)s(0)v_2 = (\alpha + \beta)Y(v_1,x)v_2.
\end{align}

\renewcommand{\theequation}{\thesection.\arabic{equation}}
\renewcommand{\thethm}{\thesection.\arabic{thm}}
\setcounter{equation}{0}
\setcounter{thm}{0}
\section{Generalized twisted modules for affine vertex operator 
algebras associated to semisimple automorphisms}

In this section, let $V$ be $M(\ell,0)$ or $L(\ell,0)$ for $\ell\ne -h^{\vee}$
(in fact, $V$ can be any quotient vertex operator algebra of $M(\ell,0)$).
Recall the automorphisms $g_{\sigma}=\mu e^{2\pi i {\rm ad}_{h}}
e^{2\pi i{\rm ad}_{\tau_{\sigma}^{-1}a_{\mathcal{N}_{g}}}}$ and 
$g_{\sigma, s}=\mu e^{2\pi i {\rm ad}_{h}}$ of 
$V$ associated to 
an automorphism $g$ of $\mathfrak{g}$ in the preceding section. 
These automorphisms of $\mathfrak{g}$ give automorphisms 
of $V$ and are still denoted using the same notations 
$g_{\sigma}$ and $g_{\sigma, s}$. 

From Proposition 3.2 in \cite{H-twisted-int}, 
we have an invertible functor $\phi_{\tau_{\sigma}}$ from the 
category of (lower-bounded or grading-restricted) generalized 
$g_{\sigma}$-twisted $V$-modules
to the category of (lower-bounded or grading-restricted) generalized  $g$-twisted $V$-modules. Thus these 
two categories are isomorphic (stronger than equivalence since the underlying vector spaces
of $W$ and $\phi_{\tau_{\sigma}}(W)$ are the same). 
In particular,  to construct (lower-bounded or grading-restricted) 
generalized  $g$-twisted $V$-modules,
we need only construct (lower-bounded or grading-restricted) generalized 
$g_{\sigma}$-twisted $V$-modules and then apply the 
functor $\phi_{\tau_{\sigma}}$.

To construct explicitly generalized 
$g_{\sigma}$-twisted $V$-modules, in this section, 
we first construct explicitly $\C/\Z$-graded 
generalized $g_{\sigma, s}$-twisted $V$-modules
from $\C/\Z$-graded  generalized $\mu$-twisted $V$-modules. 
Note that many $\C/\Z$-graded generalized $\mu$-twisted $V$-modules have been 
constructed and studied extensively since they are in fact modules 
for the corresponding twisted affine Lie algebras associated to 
the diagram automorphism $\mu$ of $\mathfrak{g}$. 
Then we discuss the construction of $\C$-graded
generalized $g_s$-twisted $V$-modules 
($g_s = e^{2 \pi i s(0)}$, see Section 5)
from generalized (untwisted) $V$-modules. 
An explicit construction of generalized 
$g_{\sigma}$-twisted $V$-modules will be given in 
the next section using the explicit construction in this section.

Let $(W, Y_{W}^{\mu})$ 
be a
$\C/\Z$-graded generalized $\mu$-twisted $V$-module
. 
We assume that $h_{W}(0)=\res_{x}x^{-1}
Y_{W}^{\mu}(h(-1)\one, x)$ acts on $W$ semisimply because $g_{\sigma, s}$ acts on $L(\ell, 0)$ 
semisimply and we always require that the semisimple part of an automorphism 
of a vertex operator algebra which acts on a twisted module is also the semisimple part of the 
action of the automorphism. Our construction in Theorem \ref{main-thm} is more general and only requires that $W$ can be decomposed as a direct sum of generalized
eigenspace of $h_W(0)$.

Let the order of $\mu$ be $r$.
 Then eigenvalues of $\mu$ are of the form $e^{2\pi i \frac{q}{r}}$ 
for $q=0, \dots, r-1$. Since $W$ is $\C/\Z$--graded,  
$W=\coprod_{n\in \C}\coprod_{q=0}^{r-1}W_{[ n]}^{[ \frac{q}{r}]},$
where $W_{[ n]}^{[\frac{q}{r}]}$ is the homogeneous subspace of $W$ 
with eigenvalue of $\mu$ being $e^{2\pi i \frac{q}{r}}$ and of conformal weight $n$. 
We define a new vertex operator map
\begin{align*}
Y_{W}^{g_{\sigma, s}}: V\otimes W&\to  W\{x\}[\log x]\nn
v\otimes w&\mapsto Y_{W}^{g_{\sigma, s}}(v, x)w
\end{align*}
by 
\begin{equation}\label{y-g-s-from-y-mu}
Y_{W}^{g_{\sigma, s}}(v, x)w=Y_{W}^{\mu}(\Delta_{V}^{(h)}(x)v, x).
\end{equation}

The eigenvalues of  $e^{2\pi i h_{W}(0)}$ are of the form 
$e^{2\pi i \alpha}$, where $\alpha\in \C$ such that
$\Re(\alpha)\in [0, 1)$. We shall denote the set of such $\alpha$ by $P_{W}^{h}$. 
Also an eigenvector of $e^{2\pi i h_{W}(0)}$ with eigenvalue $e^{2\pi i \alpha}$ 
is also an eigenvector of $h_{W}(0)$ with eigenvalue $\alpha+p$ for some $p\in \Z$. 
Since $h$ is fixed by $\mu$, $h_{W}(0)$ commutes with $\mu$.
Therefore for $q=0, \dots, r-1$ and $n\in \C$,  $W_{[ n]}^{[\frac{q}{r}]}$ 
 can be further decomposed into direct sums of eigenspaces 
 $W_{[n]}^{[ \frac{q}{r}],[ m]}$  with eigenvalue $m\in P_{W}^{h}+\Z$.  Then we have
$$W=\coprod_{n\in \C}\coprod_{q=0}^{r-1}\coprod_{m\in P_{W}^{h}+\Z}
W_{[ n]}^{[ \frac{q}{r}],[ m]}.$$

By the proof of Theorem \ref{main-thm-1} we also have  $Y_{W}^{g_{\sigma, s}}(\omega, x)\in {\rm End}\;W[[x, x^{-1}]]$
and if we write 
$$Y_{W}^{g_{\sigma, s}}(\omega, x)=\sum_{n\in \Z}L_{W}^{g}(n)x^{-n-2},$$
then by (\ref{new-L(0)}) we have 
$$L_{W}^{g_{\sigma, s}}(0)=L^{\mu}_{W}(0)
-h_{W}(0)+\frac{1}{2}(h, h)\ell,$$
where $L_{W}^{\mu}(0)$ is the corresponding Virasoro operator for the 
grading-restricted 
generalized $\mu$-twisted $L(\ell, 0)$-module $W$.
Thus we have that nonzero elements of $W^{[\frac{q}{r}],[ m]}_{[n-m+\frac{1}{2}(h, h)\ell]}$
are generalized eigenvectors of $L_{W}^{g_{\sigma, s}}(0)$ with eigenvalue $n$ 
and are also eigenvectors of $g_{\sigma, s}=\mu e^{2\pi i h_{W}(0)}$ with eigenvalue 
$e^{2\pi i (\frac{q}{r}+m)}$. Let 
$$P_{W}^{g_{\sigma, s}}=\{\beta\in \C\;|\;\Re(\beta)\in [0, 1),\; 
\beta=\frac{q}{r}+\alpha \mod \Z,\; {\rm for} \;q\in \{0, \dots, r-1\}, \; \alpha\in P_{W}^{h}\}.$$
Let 
$$W_{\langle n\rangle }^{\langle \beta \rangle}=\coprod_{\substack{q\in \{0, \dots, r-1\},\; m\in P_{W}^{h}+\Z\\
\frac{q}{r}+m=\beta{\scriptsize \mod} \Z}}W^{[ \frac{q}{r}],[ m]}_{[n+m-\frac{1}{2}(h, h)\ell]}$$
for $n\in \C$ and $\beta\in P_{W}^{g_{\sigma, s}}$. Then we have that
the space $W$ is doubly graded by the eigenvalues of
 $L_{W}^{g_{\sigma, s}}(0)$ and $g_{\sigma, s}$, that is, 
$W_{\langle n \rangle}^{\langle\beta\rangle}$ 
is the intersection of the generalized eigenspace
of  $L_{W}^{g_{\sigma, s}}(0)$ with the eigenvalues $n$ and 
the eigenspace of $g_{\sigma, s}$ with eigenvalue $e^{2\pi i \beta}$ and
\begin{equation}\label{semisimple-grading}
W=\coprod_{n\in \C}\coprod_{\beta\in P_{W}^{g_{\sigma, s}}}
W_{\langle n\rangle}^{\langle\beta\rangle}.
\end{equation}

Conversely, let $(W, Y^{g_{\sigma, s}}_{W})$
be a generalized $g_{\sigma, s}$-twisted 
$V$-module. We define
\begin{equation}\label{y-mu-from-y-g-s}
Y_{W}^{\mu}(v, x)=Y_{W}^{g}(\Delta_{V}^{h}(x)^{-1}v, x).
\end{equation}
In this case, 
$$W=\coprod_{n\in \C}\coprod_{\beta\in P_{W}^{g_{\sigma, s}}}
W_{[n]}^{[\beta]}.$$
Since $\mu$ and $g_{\sigma, s}$ commute,
$W_{[n]}^{[\beta]}$ can be decomposed into a direct sum of
eigenspaces $W_{[n]}^{[\beta], [\frac{q}{r}]}$ of $\mu$ for $q=0, \dots, r-1$
and we have 
$$W=\coprod_{n\in \C}\coprod_{q=0}^{r-1}
\coprod_{\beta\in P_{W}^{g_{\sigma, s}}}
W_{[n]}^{[\beta], [\frac{q}{r}]}.$$
Let 
$$W_{\langle n\rangle}^{\langle \frac{q}{r}\rangle}
=\coprod_{\beta\in P_{W}^{g_{\sigma, s}}}
W_{[n-m+\frac{1}{2}(h, h)\ell]}^{[\beta], [\frac{q}{r}]}$$
for $n\in \C$ and $q=0, \dots, r-1$. 
Then 
\begin{equation}\label{mu-grading}
W=\coprod_{n\in \C}\coprod_{q=0}^{r-1}
W_{\langle n\rangle}^{\langle \frac{q}{r}\rangle}.
\end{equation}

Now, applying Theorem \ref{main-thm-1} to the discussion above, we have the following result:

\begin{thm}
Let 
$\left(W, Y_{W}^{\mu}\right)$
be a $\C/\Z$-graded generalized $\mu$-twisted $V$-module. 
Assume that $h_{W}(0)=\res_{x}x^{-1}
Y_{W}^{\mu}(h(-1)\one, x)$ acts on $W$ semisimply.
The pair $(W, Y_{W}^{g_{\sigma,s}})$, with $W$ equipped with the 
new double gradings given by (\ref{semisimple-grading}) and $Y_{W}^{g_{\sigma,s}}$ defined by
(\ref{y-g-s-from-y-mu}), 
is a $\C/\Z$-graded generalized
$g_{\sigma,s}$-twisted $V$-module. Conversely,
let $\left(W, Y^{g_{\sigma, s}}_{W}\right)$
be a $\C/\Z$-graded generalized $g_{\sigma, s}$-twisted 
$V$-module. 
Then the pair $(W, Y_{W}^{\mu})$, with $W$ equipped with the 
new double gradings given by (\ref{mu-grading}) and $Y_{W}^{g_{\sigma,s}}$ defined by
(\ref{y-mu-from-y-g-s}), 
is a $\C/\Z$-graded generalized
$\mu$-twisted $V$-module.
\end{thm}

We now give explicit examples of the vertex operators $Y_{W}^{g_{\sigma,s}}(v,x)$ for $v$ in a set of generators of $V$.   
Let $a_1,\dots, a_n \in \g$ be eigenvectors of ad$_h$ with eigenvalues $\lambda_1,\dots,\lambda_n \in \C$, respectively. That is, we have
$$
\mathrm{ad}_h(a_i) = [h,a_i]= \lambda_i a_i
$$
for $1 \le i \le n$. We note here that if $\lambda_i \neq 0$, we have that 
\begin{align}
(h,a_i) = \frac{1}{\lambda_i} (h,[h,a_i]) = \frac{1}{\lambda_i} ([h,h],a_i) = 0.
\end{align}
Hence, $(h,a_i) \neq 0$ is only possible when $\lambda_i=0$. For any element $b \in \g$ with $\mu(b) = e^{\frac{2 \pi i j}{k}}b$ we will write
\begin{equation}
Y^{\mu}_W(b(-1)\one,x) = \sum_{n \in \frac{j}{k}+\Z} b_W^{\mu}(n)x^{-n-1}
\end{equation}
and
\begin{equation}\label{ss-modes}
Y^{g_{\sigma_s}}_W(b(-1)\one,x) = \sum_{n \in \C} b_W^{g_s}(n)x^{-n-1}.
\end{equation}
We note that since $\mu(h) = h$, $b$ can be decomposed into a sum of eigenvectors of $ad_h$, so without loss of generality we assume that $b$ is an eigenvector of $ad_h$ with eigenvalue $\lambda \in \{\lambda_1,\dots , \lambda_n\}$.
 We begin by explicitly computing $\Delta^{(h)}_V(x)$. We have that 
\begin{align}
-\int_{0}^{-x}Y^{\le -2}(h(-1)\one,x) &= -\int_{0}^{-x} \sum_{n \in \Z_+} h(n)y^{-n-1}\\
&= \sum_{n \in \Z_+} \frac{h(n)}{n} (-x)^{-n}
\end{align}
and so we have that 
\begin{align}
\Delta^{(h)}_V(x) = x^{-h(0)} \mathrm{exp}\left(\sum_{n \in \Z_+} \frac{h(n)}{n} (-x)^{-n}\right).
\end{align}
For each element $a_i(-1)\one$ with $1 \le i \le n$, we have that
\begin{align}
\Delta^{(h)}_V(x)b(-1)\one &= x^{-h(0)}\mathrm{exp}\left(\sum_{n \in \Z_+} \frac{h(n)}{n} (-x)^{-n}\right)b(-1)\one\\
&=  x^{-\lambda}\left(1 + h(1)(-x)^{-1}\right)b(-1)\one\\
&= b(-1)\one x^{-\lambda} - (h,b)\ell\one x^{-\lambda-1} 
\end{align}
From this, we immediately obtain 
\begin{align}
Y^{g_{\sigma_s}}_W(b(-1)\one,x) = Y^\mu_W(\Delta^{(h)}_V(x)b(-1)\one,x) = Y^\mu_W(b(-1)\one,x)x^{-\lambda} - (h,b) \ell x^{-\lambda-1}
\end{align}
In particular, since $(h,b) = 0$ if $\lambda \neq 0$, we have 
\begin{align}
Y^{g_{\sigma_s}}_W(b(-1)\one,x) = 
\begin{cases}
Y^\mu_W(b(-1)\one,x)x^{-\lambda} \  \mathrm{ if }\  \lambda \neq 0\\
Y^\mu_W(b(-1)\one,x) - (h,b)\ell x^{-1} \  \mathrm{ if }\  \lambda = 0.
\end{cases}
\end{align}
Thus we have that, for $n \in \lambda + \frac{j}{k} + \Z$, 
\begin{align}
b^{g_{\sigma_s}}_W(n) = 
\begin{cases}
b^{\mu}_W(n-\lambda) \  \mathrm{ if }\  \lambda \neq 0\\
b^\mu_W(n) \  \mathrm{ if }\  \lambda= 0 \ \mathrm{and} \ n \neq 0\\
b^\mu_W(0) -(h,b) \ell \ \mathrm{ if } \ \lambda = 0 \ \mathrm{and} \ n = 0
\end{cases}
\end{align}
In particular, we have $Y^{g_{\sigma_s}}_W(b(-1)\one,x) = \sum_{n \in \frac{j}{k}+ \lambda + \Z} b^{g_{\sigma_s}}_W(n) x^{-n-1}$. 

Instead of direct computation using the definition of $\Delta_V^{(h)}(x)$, we 
 use the $L(-1)$-derivative property to compute 
 $Y^{g_{\sigma_s}}_W(b(-n-1)\one,x)$ for $n \in \Z$, $n > 1$. 
 We recall that  
\begin{equation}
[L(m),b(n)] = -nb(m+n)
\end{equation}
for all $b \in \g$ and that $L(-1)\one = 0$ (cf. \cite{LL})). From this, it immediately follows that
\begin{equation}
b(-n-1)\one = \frac{1}{n!}{L(-1)^n}b(-1)\one
\end{equation}
We now have
\begin{align}
Y^{g_{\sigma_s}}_W(b(-n-1)\one,x) &= \frac{1}{n!} \frac{\partial^n}{\partial x^n}Y^{g_{\sigma_s}}_W(b(-1)\one,x)\\
&= \begin{cases}
\frac{1}{n!} \left(\frac{\partial}{\partial x}\right)^n\left(Y^\mu_W(b(-1)\one,x)x^{-\lambda}\right) \  \mathrm{ if }\  \lambda \neq 0\\
\frac{1}{n!} \left(\frac{\partial}{\partial x}\right)^n\left(Y^\mu_W(b(-1)\one,x) - (h,b)\ell x^{-1}\right) \  \mathrm{ if }\  \lambda = 0\\
\end{cases}
\end{align}
In particular, if $\lambda \neq 0$ we have
\begin{equation}
Y^{g_{\sigma_s}}_W(b(-n-1)\one,x) = \frac{1}{n!}\sum_{k=0}^n {-\lambda \choose k} k! \left(\left( \frac{\partial}{\partial x} \right)^{n-k} Y^\mu_W(b(-1)\one,x)\right) x^{-\lambda - k}
\end{equation}
and if $\lambda = 0$ we have:
\begin{equation}
Y^{g_{\sigma_s}}_W(b(-n-1)\one,x) = \left(\frac{1}{n!} \left( \frac{\partial}{\partial x} \right)^nY^\mu_W(b(-1)\one,x)\right) - (-1)^n (h,b) \ell x^{-n-1}
\end{equation}

Next, we discuss the construction of $\C$-graded 
generalized $g_s$-twisted modules
(recalling $g_s = e^{2 \pi i s(0)}$ from the preceding section) 
from generalized (untwisted) $V$-modules. 

Let $(W,Y_W)$ be a  generalized $V$-module. Then 
$$
Y_W(b(-1){\bf 1},x) = \sum_{n \in \Z} b_W(n)x^{-n-1}
$$
for an arbitrary $b \in \g$, where $b_{W}(n)$ for $n\in \Z$ are actions
of $b(n)$ on $W$. We assume 
that $s_W(0)$ acts on $W$ semisimply. We note here that, $W$ is a 
$1_{V}$-twisted module and the additional grading given by $1_{V}$ is trivial. 
That is, we have that $W^{[0]} = W$ and $W^{[0]}_{[n]} = W_{[n]}$, and in particular, $W$ is a $\C/\Z$-graded 
generalized $1_{V}$-twisted $V$-module.
In this case, $\mathcal{S}_g = \widetilde{\mathcal{S}}_g 
=\mathcal{N}_g=0$ since $g = e^{2 \pi i (\mathcal{S}_g+\mathcal{N}_g)} = 1_V$. 

Let 
\begin{equation}\label{y-s-from-y-1}
Y_W^{g_s}(v, x) = Y_W(\Delta_V^{(s)}(x)v, x)
\end{equation}
for $v\in V$. 
We also need an additional $\C$-grading to check that the conditions of Theorems \ref{C-graded} and \ref{strongly-C-graded} are satisfied. 
In this case, since $L_W(0)$ and $s_W(0)$ commute on $W$, we have
$$
W =W^{[0]}=\coprod_{n \in \C}
\coprod_{\alpha \in P^{g_s}_W + \mathbb{Z}} 
W^{[0], [\alpha]}_{[n]}
$$
where 
$$
W^{[0],[\alpha]}_{[n]} = \{w \in W \mid L_W(0)w = nw, \ s_{W}(0)w 
= \alpha w\}.
$$
The we have an additional $\C$-grading on $W$ given by
$$
W= \coprod_{\alpha \in P^{g_s}_W + \mathbb{Z}} W^{[0],[\alpha]},
$$
where 
$$W^{[0],[\alpha]}=\coprod_{n \in \C}
W^{[0], [\alpha]}_{[n]}.
$$
We define
$$
W_{\langle n \rangle}^{\langle \alpha \rangle} = W^{[0],[\alpha]}_{[n+\alpha - \frac{1}{2}(s,s)\ell]}
$$
so that
\begin{equation}\label{example-ss-grading}
 W  = \coprod_{n \in \C, \alpha \in P^{g_s}_W + \mathbb{Z}} W_{\langle n \rangle}^{\langle \alpha \rangle} 
\end{equation}

By Theorems \ref{C-graded} and \ref{strongly-C-graded}, we have:

\begin{thm} 
Let $(W,Y_W)$ be a generalized 
$V$-module. Assume 
that $s_W(0)$ acts on $W$ semisimply. 
Then the pair $(W, Y_{W}^{g_{s}})$, with $W$ equipped with the 
new double gradings given by (\ref{example-ss-grading}) and $Y_{W}^{g_{s}}$ defined by
(\ref{y-s-from-y-1}), is
a strongly $\C$-graded generalized $g_s$-twisted $V$-module. 
If $(W,Y_W)$  is grading restricted, then  $(W, Y_{W}^{g_{s}})$
is strongly $\C$-graded (grading restricted). 
\end{thm}
\pf
We need only check that the conditions in Theorems \ref{main-thm}, 
\ref{C-graded} and \ref{strongly-C-graded} are satisfied. 
In this case, we have $g = 1_V$ and 
$h = s(-1){\bf 1}$ in Theorem \ref{main-thm}. So it is clear that $g(h)=h$ i.e. $1_V(s(-1){\bf 1}) = s(-1){\bf 1}$. 
Also,  $s(-1){\bf 1} \in V_{(1)}$ and we have, by the commutator formula, that
\begin{equation*}
L(1)s(-1){\bf 1} = s(-1)L(1){\bf 1} + s(0){\bf 1} = 0.
\end{equation*}
Since $\widetilde{\mathcal{S}}_g = 0$, the conditions of Theorem \ref{C-graded} are trivially met. Moreover, using the notation $Q$ and $Q_n^\alpha$ from Theorem \ref{strongly-C-graded}, we have that $Q = \{ 0 \}$, and thus is a finite set, so that 
$(W, Y_{W}^{g_{s}})$ is strongly $\C$-graded.
\epfv

\begin{exam}
{\rm If $V = L(\ell,0)$, let $W = L(\ell,\lambda)$ be the irreducible quotient of 
$$\text{Ind}_{\g}^{\hat{\g}}(L(\lambda)) = U(\widehat{\g})\otimes_{U(\g \oplus \C{\bf k} \oplus \widehat{\g}_+)}L(\lambda),$$ where $L(\lambda)$ is an irreducible $\g$-module with highest weight $\lambda$ and ${\bf k}$ and $\widehat{\g}_+$ act trivially. In fact, $s(0)$ has a natural action on $L(\lambda)$ given by the action of $s$,  and thus $L(\lambda)$ can be decomposed into $s$-eigenspaces. Moreover, the maximal proper submodule $J(\ell,\lambda)$ of $M(\ell,\lambda)$ is a $\hat{\g}$-submodule of $\text{Ind}_{\g}^{\hat{\g}}(L(\lambda))$, and thus is preserved by $s(0)$. Hence, $L(\ell,\lambda)$ can be decomposed into $s(0)$-eigenspaces. In fact, since $[L_W(0),s(0)]=0$, we have that the $s(0)$-grading is compatible with the conformal weight grading, and we thus have
\begin{align}
L(\ell,\lambda) = \coprod_{m \in \C, \alpha \in \C} L(\ell,\lambda)_{\langle m\rangle}^{\langle \alpha \rangle}.
\end{align}}
\end{exam}

We conclude this section with an example to show that 
for $s\in \g$ such that $\text{ad}_{s}$ is semisimple on $\g$,
$\Delta^{s(-1)\one}_W(x)$ 
need not map a grading-restricted twisted module to another grading-restricted twisted module.

\begin{exam}
{\rm Let $\g$ be a finite dimensional simple Lie algebra, and consider an $\mathfrak{sl}(2)$-triple $\{e_\alpha, f_\alpha, h_\alpha\}$ where
\begin{align*}
[e_\alpha,f_\alpha] = h_\alpha, \ \ [h_\alpha,e_\alpha] = 2 e_\alpha, \ \ [h_\alpha,f_\alpha] = -2 f_\alpha
\end{align*}
Let $s = \frac{1}{2}h_\alpha$ so that $[s,e_\alpha] = e_\alpha$ and consider the automorphism $g_s = e^{2 \pi i s(0)}$ of $V = M(\ell,0)$. Let $W = M(\ell,0)$. We note that $W$ trivially satisfies the grading-restriction condition. We have that 
\begin{equation*}
L^{g_s}_W(0) = L_W(0) - s(0) + \frac{1}{2}(s,s) \ell.
\end{equation*}
We note that $s(0) e_\alpha(-1)^k {\bf 1} = k e_{\alpha}(-1)^k {\bf 1}$ and that $L(0) e_{\alpha}(-1)^k{\bf 1}= ke_\alpha(-1)^k{\bf 1}$, and so 
\begin{align*}
L^{g_s}_W(0) x_\alpha(-1)^k {\bf 1} = \frac{1}{2}(s,s)\ell x_\alpha(-1)^k {\bf 1}
\end{align*}
for all $k \ge 0$. Thus, we have that $\{ e_\alpha(-1)^k {\bf 1} \vert k \ge 0 \} $ is a linear independent subset of $W_{\langle \frac{1}{2}(s,s)\ell \rangle}$ and so the module $(Y^{g_s}_W,W)$ does not satisfy the grading restriction condition.}
\end{exam}

\renewcommand{\theequation}{\thesection.\arabic{equation}}
\renewcommand{\thethm}{\thesection.\arabic{thm}}
\setcounter{equation}{0}
\setcounter{thm}{0}
\section{Twisted modues for affine vertex operator algebras associated to general
 automorphisms}

In this section, $V$ is still $M(\ell,0)$ or $L(\ell, 0)$ (or a 
quotient vertex algebra of $M(\ell,0)$) for $\ell\ne -h^{\vee}$ as in the preceding section. 
We use the automorphisms of subcategories of 
generalized twisted modules  in Section 4 to 
obtain $\C/\Z$-graded 
generalized  $g_{\sigma}$-twisted $V$-modules
from generalized $g_{\sigma, s}$-twisted $V$-modules 
in this section (see Section 5 for the automorphisms $g_{\sigma}$ and 
$g_{\sigma, s}$ of $V$). 
We then discuss the construction of $\C$-graded generalized 
$g_a$-twisted modules, where $g_a$ is an inner 
automorphism $g_a = e^{2 \pi i a(0)}$ (see Section 5) which in general might 
not be semisimple.

Let $(W,Y^{g_{\sigma,s}}_W)$ be a $\C/\Z$-graded 
generalized $g_{\sigma,s}$-twisted $V$-module with
$$ W = \coprod_{n \in \C}\coprod_{\gamma \in P_{W}^{\g_{\sigma, s}}}W^{[\gamma]}_{[n]}.$$ 
Assume that $g_{\sigma,s}$ acts on $W$ semisimply and 
$(\tau_\sigma^{-1}a_{\mathcal{N}_g})_{W}(0)$ on $W$ is locally nilpotent. 
Since $\tau_\sigma^{-1}a_{\mathcal{N}_g}\in \g$, 
$(\tau_\sigma^{-1}a_{\mathcal{N}_g})(-1)\one\in V_{(1)}$. 
Let $u=(\tau_\sigma^{-1}a_{\mathcal{N}_g})(-1)\one$. 
We define a new vertex operator map
\begin{align*}
Y_{W}^{g_{\sigma}}: V\otimes W&\to  W\{x\}[\log x]\nn
v\otimes w&\mapsto Y_{W}^{g_{\sigma}}(v, x)w
\end{align*}
by 
\begin{equation}\label{y-g-from-y-g-s}
Y_{W}^{g_{\sigma}}(v, x)w=Y_{W}^{g_{\sigma,s}}(\Delta_{V}^{(u)}(x)v, x)
\end{equation}
for $v\in V$. 

Since $g_{\sigma, s}$ and $e^{\tau_\sigma^{-1}a_{\mathcal{N}_g}}$
commute and $(\tau_\sigma^{-1}a_{\mathcal{N}_g})_{W}(0)$ is locally nilpotent, 
elements of $W^{[\gamma]}_{[n]}$ for $n\in \C$ and 
$\gamma \in P_{W}^{\g_{\sigma, s}}$ are also generalized 
eigenvectors of $(\tau_\sigma^{-1}a_{\mathcal{N}_g})_{W}(0)$ with 
eigenvalue $0$. Then $P_{W}^{e^{\tau_\sigma^{-1}a_{\mathcal{N}_g}}}
=\{0\}$ and 
$$
W^{[ \gamma ], [0]}_{[ n ]}=W^{[\gamma]}_{[ n ]}
$$
for
$\gamma \in P_{W}^{\g_{\sigma, s}}$ and $n\in \C$. 
From  (\ref{new-L(0)}), we have
$$
L^{g_\sigma}_W(0) = L^{g_{\sigma,s}}_W(0) - (\tau_\sigma^{-1} 
a_{\mathcal{N}_g})(0) + \frac{1}{2}(\tau_\sigma^{-1} a_{\mathcal{N}_g},\tau_\sigma^{-1} a_{\mathcal{N}_g})\ell.
$$
Since $(\tau_\sigma^{-1} a_{\mathcal{N}_g})_{W}(0)$ is locally nilpotent,
elements of $W^{[ \gamma ]}_{[ n ]}$ are generalized eigenvectors 
of $g_\sigma$ with eigenvalue $\gamma$ and generalized eigenvectors of $L^{g_\sigma}_W(0)$ with eigenvalue 
$n + \frac{1}{2}(\tau_\sigma^{-1} a_{\mathcal{N}_g},\tau_\sigma^{-1} a_{\mathcal{N}_g})\ell$. For  $\gamma \in P_{W}^{\g_{\sigma, s}}$, 
we thus define
\begin{equation*}
W^{\langle\gamma\rangle}_{\langle n \rangle} = W^{[\gamma]}_{[n -  \frac{1}{2}(\tau_\sigma^{-1} a_{\mathcal{N}_g},\tau_\sigma^{-1} a_{\mathcal{N}_g})\ell]}
\end{equation*}
so that 
\begin{equation}\label{Lie-cz-grading}
W = \coprod_{n \in \C}\coprod_{\gamma \in P_{W}^{\g_{\sigma, s}}}
W^{\langle \gamma\rangle}_{\langle n \rangle}.
\end{equation}

Now, as a application of Theorem \ref{main-thm-1}, we obtain:

\begin{thm}
Let $(W, Y_{W}^{g_{\sigma,s}})$ be a $\C/\Z$-graded
generalized $g_{\sigma,s}$-twisted $V$-module. 
Assume that $g_{\sigma,s}$ acts on $W$ semisimply and 
$(\tau_\sigma^{-1}a_{\mathcal{N}_g})_{W}(0)$ on $W$ is locally nilpotent. Then
the pair $(W, Y_{W}^{g_{\sigma}})$, with the grading of $W$ given
 (\ref{Lie-cz-grading}) and with $Y_{W}^{g_{\sigma}}$ defined by
(\ref{y-g-from-y-g-s}), 
is a $\C/\Z$-graded generalized
$g_{\sigma}$-twisted $V$-module. 
\end{thm}
\pf
To apply Theorem \ref{main-thm-1}, it only remains to check that  $g_{\sigma,s}(\tau_{\sigma}^{-1} a_{\mathcal{N}_g}) = \tau_{\sigma}^{-1} a_{\mathcal{N}_g}$.
Indeed, we have 
\begin{align*}
g_\sigma (\tau_{\sigma}^{-1} a_{\mathcal{N}_g}) &= \tau_\sigma^{-1} g \tau_\sigma  (\tau_{\sigma}^{-1} a_{\mathcal{N}_g})\\
&= \tau_{\sigma}^{-1} a_{\mathcal{N}_g}.
\end{align*}
\epfv





We now explore the structure of the twisted module $(W,Y^{g_\sigma}_W)$. As in the previous section, consider an element $b \in \g$ which is an eigenvector of $\mu$ with eigenvalue $e^{2 \pi i \frac{j}{k}}$ and an eigenvector of $\text{ad}_h$ with eigenvalue $\lambda$. Importantly, we have that $b$ is a generalized eigenvector of $g_\sigma$ and $e^{2 \pi i (\tau_\sigma^{-1} a_{\mathcal{N}_g})(0)}$.
 Since $b$ is a generalized eigenvector of $g_{\sigma}$, 
 $\text{ad}_{\tau_\sigma^{-1} a_{\mathcal{N}_g}}$ on $b$ is nilpotent.
Then there exists $M \in \mathbb{Z}_+$ such that $(\text{ad}_{\tau_\sigma^{-1}
 a_{\mathcal{N}_g}})^M (b) = 0$. We write
$$
Y^{g_{\sigma}}_W(b(-1)\one,x) =\sum_{k \ge 0} \sum_{n \in \C} b^{g_{\sigma}}(n,k)x^{-n-1} (\log x)^k
$$

We begin by explicitly computing 
$\Delta^{(u)}_V(x)=\Delta^{(\tau_\sigma^{-1} a_{\mathcal{N}_g})(0))}_V(x)$. 
As above, we have
$$
-\int_{0}^{-x}Y^{\le -2}((\tau_\sigma^{-1} a_{\mathcal{N}_g})(-1)\one,x) = \sum_{m \in \Z_+} \frac{(\tau_\sigma^{-1} a_{\mathcal{N}_g})(m)}{m} (-x)^{-m}
$$
and so we have that 
$$
\Delta^{(\tau_\sigma^{-1} a_{\mathcal{N}_g})}_V(x) = \exp(-(\tau_\sigma^{-1} a_{\mathcal{N}_g})(0) \log x) \mathrm{exp}\left(\sum_{m \in \Z_+} \frac{(\tau_\sigma^{-1} a_{\mathcal{N}_g})(m)}{m} (-x)^{-m}\right).
$$
Thus we have that 
\begin{align*}
\Delta&^{((\tau_\sigma^{-1} a_{\mathcal{N}_g})(0))}_V(x) b(-1)\one \nn
&= \exp(-(\tau_\sigma^{-1} a_{\mathcal{N}_g})(0) \log x) \mathrm{exp}\left(\sum_{m \in \Z_+} \frac{(\tau_\sigma^{-1} a_{\mathcal{N}_g})(m)}{m} (-x)^{-m}\right) b(-1)\one \\
&= \exp(-(\tau_\sigma^{-1} a_{\mathcal{N}_g})(0) \log x) \left(b(-1)\one - ((\tau_\sigma^{-1} a_{\mathcal{N}_g}),b)\ell \one x^{-1}\right)\\
&= \sum_{l\ge0} \frac{1}{l!}(-(\tau_\sigma^{-1} a_{\mathcal{N}_g})(0)\log x)^l  b(-1)\one - ((\tau_\sigma^{-1} a_{\mathcal{N}_g}),b)\ell \one x^{-1}\\
&= \sum_{l=0}^{M} \frac{1}{l!}(-1)^l (\text{ad}_{(\tau_\sigma^{-1} a_{\mathcal{N}_g})}^l(b))(-1)\one (\log x)^l - ((\tau_\sigma^{-1} a_{\mathcal{N}_g}),b)\ell \one x^{-1}.
\end{align*}
We write 
$$
Y^{g_{\sigma_s}}_W((\text{ad}_{(\tau_\sigma^{-1} a_{\mathcal{N}_g})}^l(b))(-1)
\one,x) =\sum_{m \in \lambda + \frac{j}{k}+ \Z} (\text{ad}^l_{(\tau_\sigma^{-1} 
a_{\mathcal{N}_g})}(b))_W^{g_{\sigma_s}}(m)x^{-m-1}.
$$
Note that there is no logarithm of $x$ in 
$Y^{g_{\sigma_s}}_W((\text{ad}_{(\tau_\sigma^{-1} 
a_{\mathcal{N}_g})}^l(b))(-1)
\one,x)$ since the automorphsim $g_{\sigma_s}$ of $V$ is semisimple. 
Thus, we have that
\begin{align*}
Y^{g_\sigma}_W(b(-1)\one,x) &= Y^{g_{\sigma_s}}_W(\Delta^{(\tau_\sigma^{-1} a_{\mathcal{N}_g})}_V(x)b(-1)\one,x)\\
&= \sum_{l=0}^{M} \frac{(-1)^l}{l!}
 Y^{g_{\sigma_s}}_W((\text{ad}_{(\tau_\sigma^{-1} 
 a_{\mathcal{N}_g})}^l(b))(-1)\one,x)(\log x)^l - (\tau_\sigma^{-1} a_{\mathcal{N}_g},b)Y^{g_s}(\ell \one,x)x^{-1}\\
&= \sum_{l=0}^{M} \sum_{m \in \lambda + \frac{j}{k} + \Z} \frac{(-1)^l}{l!} 
(\text{ad}^n_{(\tau_\sigma^{-1} a_{\mathcal{N}_g})}(b))_W^{g_{\sigma_s}}(m) x^{-m-1} (\log x)^l - (\tau_\sigma^{-1} a_{\mathcal{N}_g},b)\ell x^{-1}
\end{align*}
In particular, we have that
\begin{align*}
b^{g_\sigma}(m,l) = \begin{cases}
{\displaystyle \frac{(-1)^l}{l!}(\text{ad}^l_{(\tau_\sigma^{-1} a_{\mathcal{N}_g})}(b))_W^{g_{\sigma_s}}(m) }
& \text{if } m \in \lambda + \frac{j}{k} + \Z \text{ and } 0 \le l \le M\\
{\displaystyle \frac{(-1)^l}{l!}(\text{ad}^l_{(\tau_\sigma^{-1} 
a_{\mathcal{N}_g})}(b))_W^{g_{\sigma_s}}(0) 
- (\tau_\sigma^{-1} a_{\mathcal{N}_g},b)\ell}
 & \text{if } m=0 \text{ and } 0 \le  l \le M\\
 0&\text{if } m\ne 0 \text{ and }m\not\in \lambda + \frac{j}{k} + \Z.
\end{cases}
\end{align*}

As in the preceding section, instead of direct computation using the definition of $\Delta_V^{(\tau_\sigma^{-1} a_{\mathcal{N}_g})}(x)$, we use the $L(-1)$-derivative property to compute $Y^{g_{\sigma}}_W(b(-n-1)\one,x)$ for $n \in \Z$, $n > 1$. Here, the calculation is straightforward, though the expression 
becomes complicated:
\begin{align*}
&Y^{g_{\sigma}}_W(b(-n-1)\one,x)\\ 
&= \frac{1}{n!} Y^{g_{\sigma}}_W(L(-1)^nb(-1)\one,x)\\
&= \frac{1}{n!}\left(\frac{\partial}{\partial x}\right)^n Y^{g_{\sigma}}_W(b(-1)\one,x)\\
&=\frac{1}{n!}\left(\frac{\partial}{\partial x}\right)^n\left(\sum_{l=0}^{M} \sum_{m \in \lambda + \frac{j}{k} + \Z} \frac{(-1)^l}{l!} (\text{ad}^l_{(\tau_\sigma^{-1} a_{\mathcal{N}_g})}(b))_W^{g_{\sigma_s}}(m) x^{-m-1} (\log x)^j - (\tau_\sigma^{-1} a_{\mathcal{N}_g},b)\ell x^{-1}\right)\\
&=\sum_{l=0}^{M} \sum_{m \in \lambda + \frac{j}{k}+ \Z}\sum_{r=0}^n  \frac{1}{n!}{n \choose r} \frac{(-1)^l}{l!} (\text{ad}^l_{(\tau_\sigma^{-1} a_{\mathcal{N}_g})}(b))_W^{g_{\sigma_s}}(m) \left( \left(\frac{\partial}{\partial x}\right)^{n-r}x^{-m-1}\right) \left(\left(\frac{\partial}{\partial x}\right)^r(\log x)^l\right)\\
& \quad  - (-1)^n(\tau_\sigma^{-1} a_{\mathcal{N}_g},b)\ell x^{-n-1}.
\end{align*}

Finally, we discuss the construction of $\C$-graded generalized 
$g_a$-twisted modules from $\C$-graded generalized $g_{s}$-twisted 
$V$-modules, where $g_a = e^{2 \pi i a(0)}$ is an inner 
automorphism of $V$ and $g_{s}= e^{2 \pi i s(0)}$ is its semisimple 
part (See Section 5). 

Suppose that 
$(W,Y_W^{g_s})$ is a $\C$-graded generalized $g_s$-twisted $V$-module. 
We assume that the actions $s_{W}(0)$ and $n_{W}(0)$ of $s(0)$ 
and $n(0)$ on $W$ are the semisimple and nilpotent parts of the action 
$a_{W}(0)$ of $a(0)$ since $s$ and $n$ are the 
semisimple and nilpotent part of the derivation $a$ of $\g$. 
In particular, 
\begin{equation}
W = \coprod_{n \in C, \alpha \in P^{g_s}_W + \Z}W_{[n]}^{[\alpha]}.
\end{equation}
We now to apply the results in Section 4 to the case that 
$V$ is $M(\ell,0)$ or $L(\ell,0)$,
$g=g_{s}$ and $u=n(-1)\one$. In this case,  $g_{u}=e^{2\pi in(0)}$ and 
$gg_{u}=e^{2\pi i a(0)}$. 

We define 
\begin{equation}\label{vo-g-a}
Y_W^{g_a}(v, x) = Y^{g_s}_W(\Delta_V^{(n(-1)\one)}(x)v,x)
\end{equation}
for $v\in V$. 

Since by assumption, $n_{W}(0)$ is locally nilpotent, 
every element of $W$ is a generalized eigenvector of 
$n_{W}(0)$ with eigenvalue $0$. Then the
$\C$-grading on $W$ given by the generalized eigenspace 
of $n_{W}(0)$ is trivial, that is, 
$$
W^{[\alpha],[0]}_{[n]} = W^{[\alpha]}_{[n]}
$$
for $n\in \C$ and $\alpha \in P_W^{g_s} + \mathbb{Z}$. 
We note also that 
$$
L^{g_a}_W(0) = L^{g_s}_W(0) - n(0) + \frac{1}{2}(n,n)\ell
$$
so that the space $W_{[n]}$ is made up of generalized eigenvectors of $L^{g_a}(0)$ with eigenvalue $n-0 + \frac{1}{2}(n,n)\ell$. We define
\begin{equation}\label{grading-g-a}
W_{\langle n\rangle}^{\langle \alpha\rangle}= 
W^{[\alpha]}_{[n-\frac{1}{2}(n,n)\ell]}.
\end{equation}
Moreover, if $W$ is grading restricted, 
we have that each $W_{\langle n\rangle}^{\langle \alpha\rangle}$ 
is clearly finite dimensional since each $W^{[\alpha]}_{[n-\frac{1}{2}(n,n)\ell]}$ is finite dimensional, and that $W_{\langle n \rangle}^{\langle \alpha\rangle} = 0$ 
when the real part of $n$ is sufficiently negative since 
$W^{[\alpha]}_{[n-\frac{1}{2}(n,n)\ell]}=0$ when the real part of $n$ 
is sufficiently negative. Thus $W$ is also grading restricted with the grading
given by (\ref{grading-g-a}). 

From Theorem \ref{main-thm-1}, we obtain the following result:

\begin{thm}
Let $(W,Y_W^{g_s})$ be a $\C$-graded generalized $g_s$-twisted $V$-module. 
such that $s_{W}(0)$ and $n_{W}(0)$ 
are the semisimple and nilpotent parts of 
$a_{W}(0)$. Then the pair $(W, Y_W^{g_a})$ with the 
grading of $W$ given by (\ref{grading-g-a}) and 
with $Y_W^{g_a}$ defined by (\ref{vo-g-a}) is a 
$\C$-graded generalized $g_a$-twisted $V$-module. 
If $(W,Y_W^{g_s})$ is grading restricted, then 
$(W, Y_W^{g_a})$ is also grading restricted.
\end{thm}
\pf
We need only check the conditions in Theorems \ref{main-thm} and 
\ref{main-thm-1}. Since $[s,n]=0$, we have that $s(0)n(-1){\bf 1} 
= n(-1)s(0){\bf 1} = 0$. Since $g_s = e^{2 \pi i s(0)}$, we have that
\begin{equation*}
g_s(n(-1){\bf 1}) = \sum_{k \ge 0}\frac{(2 \pi i s(0))^k}{k!} n(-1){\bf 1} = n(-1){\bf 1}.
\end{equation*}
By the commutator formula, we also have
\begin{equation*}
L(1)n(-1){\bf 1} = n(-1)L(1){\bf 1} + n(0){\bf 1} = 0.
\end{equation*}
Hence, all the conditions of Theorems \ref{main-thm} and \ref{main-thm-1}
are satisfied and so our theorem follows.
\epfv

\noindent {\small \sc Department of Mathematics, Rutgers University,
110 Frelinghuysen Rd., Piscataway, NJ 08854-8019}

\noindent {\em E-mail address}: yzhuang@math.rutgers.edu

\vspace{1em}

\noindent {\small \sc  Department of Mathematics and Computer Science, Ursinus College, 
Collegeville, PA 19426-1000} 

\noindent {\em Email address}: csadowski@ursinus.edu

\end{document}